\documentclass[letterpaper, 11pt]{article}
\usepackage{amsmath}
\usepackage{amssymb}
\usepackage{amsfonts}
\usepackage{mathrsfs}
\usepackage{bbm}
\usepackage{color}
\usepackage{graphicx}
\usepackage{epsfig}

\definecolor{Red}{rgb}{1,0,0}
\definecolor{Blue}{rgb}{0,0,1}
\definecolor{Olive}{rgb}{0.41,0.55,0.13}
\definecolor{Green}{rgb}{0,1,0}
\definecolor{MGreen}{rgb}{0,0.8,0}
\definecolor{DGreen}{rgb}{0,0.55,0}
\definecolor{Yellow}{rgb}{1,1,0}
\definecolor{Cyan}{rgb}{0,1,1}
\definecolor{Magenta}{rgb}{1,0,1}
\definecolor{Orange}{rgb}{1,.5,0}
\definecolor{Violet}{rgb}{.5,0,.5}
\definecolor{Purple}{rgb}{.75,0,.25}
\definecolor{Brown}{rgb}{.75,.5,.25}
\definecolor{Grey}{rgb}{.5,.5,.5}
\definecolor{Black}{rgb}{0,0,0}

\newcommand{\acal}{\mathcal{A}}
\newcommand{\bcal}{\mathcal{B}}

\newcommand{\dcal}{\mathcal{D}}
\newcommand{\ecal}{\mathcal{E}}
\newcommand{\fcal}{\mathcal{F}}
\newcommand{\gcal}{\mathcal{G}}

\newcommand{\kcal}{\mathcal{K}}
\newcommand{\lcal}{\mathcal{L}}

\newcommand{\rcal}{\mathcal{R}}
\newcommand{\scal}{\mathcal{S}}

\newcommand{\qcal}{\mathcal{Q}}

\newcommand{\real}{\mathbb{R}}

\newcommand{\nintgr}{\mathbb{N}}

\newcommand{\eps}{\varepsilon}

\newcommand{\bdm}{\begin{displaymath}}
\newcommand{\edm}{\end{displaymath}}
\newcommand{\bea}{\begin{eqnarray*}}
\newcommand{\eea}{\end{eqnarray*}}
\newcommand{\bean}{\begin{eqnarray}}
\newcommand{\eean}{\end{eqnarray}}

\newcommand{\bfx}{\mathbf{x}}
\newcommand{\bfy}{\mathbf{y}}

\newcommand{\prob}{\mathbb{P}}
\newcommand{\expec}{\mathbb{E}}
\newcommand{\var}{\mathrm{Var}}

\newcommand{\depth}{\Delta}

\newcommand{\poly}{\mathrm{poly}}

\newcommand{\diff}{\mathrm{d}}

\newcommand{\undtheta}{\underline{\theta}}
\newcommand{\ovtheta}{\overline{\theta}}

\newcommand{\delhat}{{\hat \delta}_{ab}}

\newcommand{\del}{\delta_{ab}}

\newcommand{\edgehat}[1]{\hat{\delta}_e (#1)}

\newcommand{\distdw}{\widehat W}

\newcommand{\distd}{\hat d}
\newcommand{\dist}{d}
\newcommand{\pathdepth}{\Delta}

\newcommand{\cluster}[1]{\widehat{H}_{#1}}
\newcommand{\vcluster}[1]{\widehat{V}_{#1}}
\newcommand{\ecluster}[1]{\widehat{E}_{#1}}

\newcommand{\vclusters}[2]{\hat{v}_{#1}^{(#2)}}

\newcommand{\ball}[2]{\widehat{B}_{#1}^{(#2)}}

\newcommand{\intersect}{\widehat{\Phi}}
\newcommand{\minicontractor}{\textsc{Mini Contractor}}
\newcommand{\extender}{\textsc{Extender}}

\newcommand{\minipart}[3]{\psi_{#1}(#2,#3)}
\newcommand{\fullpart}[3]{\bar\psi_{#1}(#2,#3)}
\newcommand{\minipartone}[3]{\psi^{(#2)}_{#1}(#2,#3)}
\newcommand{\miniparttwo}[3]{\psi^{(#3)}_{#1}(#2,#3)}
\newcommand{\fullpartone}[3]{\bar\psi^{(#2)}_{#1}(#2,#3)}
\newcommand{\fullparttwo}[3]{\bar\psi^{(#3)}_{#1}(#2,#3)}
\newcommand{\partsize}[2]{r(#1,#2)}

\newtheorem{theorem}{Theorem}
\newtheorem{proposition}{Proposition}

\newtheorem{definition}{Definition}

\newtheorem{lemma}{Lemma}

\newtheorem{remark}{Remark}

\newtheorem{assumption}{Assumption}

\newenvironment{proof}{\noindent{\textbf{Proof:}}}{$\blacksquare$\vskip\belowdisplayskip}

\author{{\bf Shankar Bhamidi}\\
Department of Statistics and OR\\
UNC Chapel Hill
\and
{\bf Ram Rajagopal}\\
Department of EECS\\
UC Berkeley
\and
{\bf S\'ebastien Roch}\\
Department of Mathematics\\
UC Los Angeles
}

\title{Network Delay Inference from Additive Metrics}

\begin{document}

\maketitle

\thispagestyle{empty}

\begin{abstract}
We use computational phylogenetic techniques
to solve a central problem in inferential network monitoring.
More precisely, we design a novel algorithm for multicast-based delay inference, that is, the 
problem of reconstructing delay characteristics of a network
from end-to-end delay measurements on network paths. 
Our inference algorithm is based on additive metric techniques used in phylogenetics. 
It runs in polynomial time and requires a sample
of size only $\poly(\log n)$. We also show how to recover the topology of the routing tree.
\end{abstract}


\section{Introduction}

\paragraph{Network tomography.} Inferential network monitoring---also 
known as network tomography~\cite{Vardi:96}---consists in reconstructing various
properties of large communication networks from indirect measurements in order to facilitate 
the management of these networks. Network inference can be achieved by two general
approaches. In the {\em internal} approach, one takes measurements directly at the edges and nodes of the network.
This approach suffers from several drawbacks: 
the network operator may not allow access to internal devices of the network or
may not make public measurements on them; 
the routers may not have the technological capabilities to perform
the required measurements;
direct measurements may create extra computational burden
as well as congestion in the network.
This has led some in the networking community to consider instead the {\em external} approach. In this case,
one uses so-called ``end-to-end'' measurements, e.g., measurements of delays 
or rate of packet drops between nodes in the network, 
and seeks to infer the desired network properties from them. 
This gives rise to an inverse problem similar to tomographic image reconstruction.

Our aim in this paper is to propose a novel approach to this problem.
We focus on multicast-based inference. Multicast routing consists in sending a packet from a source
to a set of receivers through a routing tree. The packet is duplicated at each branch point and sent further down the tree.
The routing tree is generally unknown to the user.
The idea is to use inherent correlation of measurements between different receivers to reconstruct the topology of the routing tree 
as well as to estimate link properties of this tree. The main link property we consider here is the delay distribution. 
The multicast inference approach was introduced in~\cite{CaDuHoTo:99, LoDuHoTo:02}.

A core difficulty of the problem is to devise efficient, 
scalable algorithms which consistently estimate the desired network properties.
Several techniques have been used in the network tomography literature, 
notably maximum pseudo-likelihood, EM algorithms and
Markov chain Monte Carlo methods. 
See~\cite{CCLNY:04} for a detailed survey and bibliographic references.
In this paper, we introduce a new methodology for multicast delay inference 
inspired by techniques from the field of {\em phylogenetics} in biology, that is, the reconstruction
of evolutionary trees from molecular data. 
Our methodology has the advantage of being provably consistent and computationally efficient. It also
uses a small asymptotic sample size. This is crucial to reduce the burden on the network
as well as to obtain a consistent ``snapshot'' of the network, which is intrinsically dynamic in nature.
Typical networks undergo sporadic medium to large-scale changes in structure over time, 
therefore algorithms with low sample complexity are essential.
Concurrently to our work, Liang et al.~\cite{LiMoYu:u} used similar ideas to tackle
the related multicast packet loss inference problem. Also, Ni and Tatikonda~\cite{NiTatikonda:06}
independently proposed a Markov-based inference algorithm similar to ours for multicast delay inference---although
our work appears to be the first rigorous analysis of the sample complexity of this approach. 
See Section~\ref{sec:results} for a precise statement of our results and Section~\ref{sec:discussion}
for a discussion of previous work.
The results detailed here were first announced in~\cite{BhRaRo:arxiv}.

\paragraph{Phylogeny background.} A core problem 
in evolutionary biology is the inference of evolutionary histories
of organisms 
from molecular data. Evolution
is usually represented by a tree where branching points indicate speciation events. 
The root of the tree is the common ancestor to all species in the tree and
the leaves are contemporary (extant) species. Molecular
data is assumed to evolve according to a standard Markov model. 
The {\it phylogenetic reconstruction} problem is the following. From measurement of sequences of molecular data
at the leaves, one seeks to reconstruct the topology
of the evolutionary tree as well as mutation characteristics along the branches. 
See~\cite{Felsenstein:04} and~\cite{SempleSteel:03} for an overview of the field of phylogenetics.

Various statistical and computational techniques 
have been used to solve the phylogenetic reconstruction problem: maximum likelihood,
bayesian, parsimony, and distance-based methods. In this paper, 
we adapt and extend {\em distance-based techniques}
to deal with a class of models introduced in~\cite{LoDuHoTo:02} 
in connection with the multicast network inference 
problem---this new class of models is similar to the Markov models 
used in phylogenetics but presents challenges of its own.
The main idea
in distance-based methods is to define a so-called tree metric 
from mutation parameters. A tree metric
is a metric on the leaves of the tree which can be realized as a path metric on a corresponding
weighted tree. (See Section~\ref{sec:phylo} for more details.) After being estimated, the metric allows the reconstruction of the tree and its 
characteristics. A main advantage of this approach is 
that it leads to computationally efficient algorithms 
with provable sample requirement guarantees.

\subsection{Basic Definitions}\label{sec:definitions}

\paragraph{A broadcasting process on a tree.} 
We now give a more formal statement of the multicast inference problem introduced in~\cite{LoDuHoTo:02}. 
Let $T=(V,E)$ be a tree on
$n+1$ leaves $L$---representing the {\em routing tree}---and let $\{d_e\}_{e\in E}$
be a set of independent positive random variables
on the edges---representing the {\em delays}. 
Leaf $0$, the {\em source}, is the root of the tree. The remaining
$n$ leaves are the {\em receivers}. We assume that all internal nodes
have degree at least 3.

A realization of the {\em multicast delay
process} works as follows: the root sends a packet to the receivers
through the routing tree; at every branching point, the packet
is duplicated; on every link $e$, 
an independent random delay $d_e$ is experienced by the packet.
More formally, we define the multicast delay process $\{D_u\}_{u\in V}$ as follows.
Let $P_{ij}$ be the path (set of edges) between
nodes $i$ and $j$ in $T$. 
For a node $u$, let
\begin{equation}\label{eq:total}
D_u = \sum_{e \in P_{0u}} d_e.
\end{equation}
Note that $D_u$ is the total delay at node $u$ in the network.

\paragraph{The multicast inference problem.}
The tree and delay distributions
are actually unknown to us. 
We are only given access to $k$ independent samples
of delays at the leaves $\{D^1_a\}_{a\in L}, \ldots, \{D^k_a\}_{a\in L}$.
Our goal is to reconstruct the routing tree and
estimate the delay distributions using these samples.
We now define more precisely what we mean by
the estimation of the delay distributions. In this work,
we assume that each edge delay distribution (in general, different) 
is characterized by a constant number,
say $J -1 > 1$ (independent of $n$),
of consecutive central moments. That is, we assume there are 
{\em characteristic moments}
\begin{equation*}
w_{e}^{(j)} = \expec\left[\left(d_e -\expec[d_e]\right)^j\right],
\end{equation*}
for all $e\in E$ and $2\leq j\leq J$.
Our goal is to estimate these moments within a fixed accuracy.
More formally, we make the following assumption. We first need a definition.
\begin{definition}[Regular Families]
Let $\eps > 0$ and $J\geq 2$ be fixed.
Let $\qcal = \{Q_\theta\}_{\theta\in \Theta}$ be a family of distributions
on $\real$ parametrized by $\theta \in \Theta$
where $\Theta$ is a subset of an Euclidean space. 
Let $\left\{w^{(j)}(\theta)\right\}_{\{2\leq j\leq J\}}$
be the first $J-1$ central moments of $Q_\theta$.
We say that the family $\qcal$
is $(\eps,J)$-regular if there exists a map $\Psi$ from $\real^{J-1}$ to $\Theta$ 
and a $\delta > 0$ such that
if the vector $\mathbf{\hat w} =  \left\{\hat w^{(j)}\right\}_{\{2\leq j\leq J\}}$ satisfies
\begin{equation*}
\left|\hat w^{(j)} - w^{(j)}(\theta)\right| \leq \delta
\end{equation*}
for all $2\leq j\leq J$, then 
\begin{equation*}
\left\| Q_\theta - Q_{\Psi(\mathbf{\hat w})}\right\|_1 \leq \eps.
\end{equation*}
\end{definition}
In Appendix~\ref{sec:regular}, we give simple examples of regular families.
\begin{assumption}[Regularity and Boundedness]\label{assump:regularity}
Let $\eps > 0$ and $J\geq 2$ be fixed (independent of $n$). 
We assume that all edge delay distributions are from a
fixed $(\eps, J)$-regular family of distributions.
Furthermore, we assume that the delays are uniformly bounded, namely 
there is a constant $M > 0$ independent
of $n$ such that for all $e\in E$, $d_e \in [0,M]$.
\end{assumption}
This framework is simple enough to be tractable yet
general enough to accommodate large classes of distributions: 
parametrized distributions, e.g., beta distributions; 
and nonparametrized
distributions, e.g., discretized distributions on $\{0,1,\ldots,M\}$.
Further we need the following assumption.
\begin{assumption}[Lower Bound on Second Moment]\label{assump:moments}
We assume that there is a constant $f > 0$ (independent of $n$) such that 
for all $e\in E$,
\begin{equation*}
w_{e}^{(2)} \geq f.
\end{equation*}
\end{assumption}

\begin{figure}[!t]
\begin{center}
\input{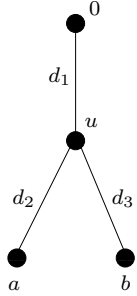}
\caption{Unidentifiability of Mean Delay: If one were to replace
$d_1$ with $d_1 + \mu$ and $d_2, d_3$ with $d_2 - \mu, d_3 - \mu$ for $\mu > 0$ 
(assuming $\mu$ can be chosen so that all delays remain positive) then
the distribution of delays at $a,b$ would be unchanged. This example also
shows that one cannot deduce the delays on all edges given total delays
at all leaves.}\label{fig:ident}
\end{center}
\end{figure}

To sum up, the {\em multicast inference problem} is defined as follows.
\begin{definition}[Multicast Inference Problem, Moment Version]
Let $\eps > 0$ and $J\geq 2$ be fixed.
The {\em multicast inference problem} consists in the following. Let $T$ and
$\left\{w_{e}^{(j)}\right\}_{\{e\in E, 2\leq j\leq J\}}$ be any tree
(with internal degrees at least 3)
and set of central moments on edges.
Given samples of delays 
at the leaves, we are required to:
\begin{enumerate}
\item {\bf Tree Reconstruction.} Recover $T$.
\item {\bf Moment Estimation.} Estimate all
characteristic moments $\left\{w_{e}^{(j)}\right\}_{\{e\in E, 2\leq j\leq J\}}$
within $\eps$.
\end{enumerate}
\end{definition}
\begin{remark}
As noted by Lo Presti et al.~\cite{LoDuHoTo:02}, the means of the
edge delay distributions are, in general, unidentifiable. See Figure~\ref{fig:ident}
for an illustration. In particular, one cannot hope to recover the {\em deterministic} 
transmission delay on each link. 
But, as noted in~\cite{LoDuHoTo:02}, this is not a major issue. Indeed, in practice, one
is only interested in the {\em variable} portion of the delay, that is, 
the portion resulting from
traffic. To restore identifiability, Lo Presti et al.~proceed by subtracting
the lowest observed delay on each receiver, in order to remove the (estimated) deterministic 
component of the delay. They further assume that the variable portion of the delay ``starts at 0.''
We also make this last assumption 
(see our examples of regular delay distributions in Appendix~\ref{sec:regular}). 
However, instead of subtracting the minimum observed delay (which may be unreliable on a large network),
we use {\em central} moments---which are not affected by the deterministic transmission delay.
\end{remark}

\subsection{Our Results.}\label{sec:results}

Our main result is the following theorem.
\begin{theorem}[Main Result]
Let $\eps > 0$ and $J \geq 2$ be fixed.
Let Assumptions~\ref{assump:regularity} and~\ref{assump:moments}
hold.
Then, there is a polynomial-time algorithm which
solves the multicast inference problem with
high probability using $k = O\left(\poly(\log n) \right)$ samples.
\end{theorem}
See Theorems~\ref{thm:rti},~\ref{thm:symer}, and~\ref{thm:generaler} below
for more precise statements. 

The proofs of the main theorems rely on the important
notion of a tree metric from phylogenetics. Roughly speaking, 
a tree metric is a metric on the leaves of a tree such that the distance
between any two leaves can be written
as a sum of edge weights on the corresponding path. (See Section~\ref{sec:phylo} for definitions.) 
There are two components to our algorithm: 
\begin{enumerate}
\item {\bf Topology reconstruction:} The reconstruction of the routing tree can be achieved by adapting 
known phylogenetic reconstruction algorithms---once the proper delay-based
metric is defined. This result is proved in Section~\ref{sec:top}. The relevant
phylogenetic background is introduced in Section~\ref{sec:phylo}. 

\item {\bf Moment estimation on edges:} Most of the technical
work of this paper is in deriving and analyzing a metric-based algorithm for inferring
edge delay distributions (Theorems~\ref{thm:symer} and~\ref{thm:generaler}). 
For this purpose, a) we relax the notion of a tree metric to allow nonnegative
edge weights, b) we define appropriate delay-based metrics, and c) we show how to
estimate these metrics. The analysis relies on large deviations arguments.
\end{enumerate}

As far as we are aware, our algorithm is the first multicast inference algorithm
to be both provably efficient and consistent. Previous work concerned mostly
non-rigorous techniques such as maximum pseudo-likelihood and EM
algorithms. See~\cite{CCLNY:04} for details. An exception is the independent,
unpublished work of Liang et al.~\cite{LiMoYu:u} which uses techniques similar
to ours in the related context of multicast packet drop inference.

\subsection{Discussion}\label{sec:discussion}

\paragraph{Validity of assumptions.} The multicast delay process defined
in Section~\ref{sec:definitions} relies on
two basic assumptions about routing and traffic which makes its analysis
possible: {\em temporal} and {\em spatial independence}. 
In reality, of course, both assumptions are violated to some extent. 
Lo Presti et al.~\cite{LoDuHoTo:02}
(see also~\cite{CaDuHoTo:99}) studied the effect of these violations
empirically and concluded that the multicast delay process is a useful first
approximation to the underlying complex process. We briefly summarize their findings.

{\em Temporal dependence}---delays at a given link being correlated at different points
in time---is common in communication networks. But, as it turns out, its impact is rather mild
for our purposes. 
Indeed the type of inference procedure studied in~\cite{LoDuHoTo:02} (as well as in the current paper) 
does not actually require independence in time
but only {\em ergodicity}---a much weaker assumption; 
more precisely, the estimator in~\cite{LoDuHoTo:02} (and in the current paper) 
is consistent as long as
the delay process is ergodic. Hence, the temporal dependencies impact only 
the {\em convergence rate} of the inference procedure. 
Lo Presti et al.~showed empirically that, although this effect cannot be ignored, it is rather mild.
Quantifying exactly the effect of temporal correlations on the theoretical convergence rate 
of an estimator is non-trivial. 

As for {\em spatial correlations}---dependencies in delays on neighboring links---Lo Presti et al.~found
that they can produce a systematic bias in the estimation. However, they showed empirically
that the bias is a small, second-order effect, possibly---they argue---because 
the diversity of traffic on the network results only in localized, short-term correlations in delays.    
They also point out that very little is known about the precise structure of such spatial
correlations in real networks, making it hard to derive a good model for them.

Another assumption implicit in our model is that the process, including the routing
tree itself, remains homogeneous over time. In fact, there are sporadic large-scale
changes in the network. These explain why
a {\em low sample complexity} is critical for an inference procedure to be useful in practice. 
Minimizing the sample complexity is the main focus of this paper.

\paragraph{Related results.}
The {\em multicast delay inference} problem was formalized by Lo Presti et al.~in~\cite{LoDuHoTo:02}.
In that paper, the authors give a procedure to infer a discretized delay distribution 
on each link, 
given the routing tree topology. Their algorithm is based on an ad-hoc fixed point equation
that is solved by least squares. Moreover, these authors show that their estimator 
is asymptotically normal with a
variance-covariance matrix depending implicitly on the delay characteristics. More explicit formulas are given
in the limit of small delays. The algorithm is tested on small networks and the dependence
on the size is not given.

More recently, Ni and Tatikonda~\cite{NiTatikonda:06,NiTatikonda:07,NiTatikonda:08,NXTY:08}---in work subsequent
to ours~\cite{BhRaRo:arxiv}---used phylogenetic techniques to recover the routing tree topology in this context.
Similarly to the current paper, they use distance-based techniques.
The basic algorithm they consider is the well-known Neighbor-Joining (NJ) algorithm which
they apply to various tree metrics, for instance, the delay variance metric (as we do here).
They also deal with trees of internal degrees higher than 3 by introducing a variant of NJ called
Rooted Neighbor-Joining (RNJ)~\cite{NiTatikonda:08} (based on a technique equivalent 
to what is known in phylogenetics as the Farris transform~\cite{Farris:73}).
They show more precisely that RNJ is a consistent estimator of the routing tree, but no convergence rate
is given. Note, however, that RNJ 
has in fact a high sample complexity
due to its reliance on the {\em diameter} of the tree. See, e.g.,~\cite{Atteson:99}.
See also our discussion about diameter v.~depth in Section~\ref{sec:distorted}.
Here, we make use of state-of-art phylogenetic reconstruction techniques to
derive a low sample complexity algorithm for routing tree reconstruction.
We also show how to infer delay distributions. A technique to infer discrete delays
was also subsequently obtained by Ni and Tatikonda~\cite{NiTatikonda:07} (although no convergence rate
is provided).

A related network tomography problem is the so-called {\em multicast link loss inference} problem,
where one observes packet losses at the receivers of a multicast routing tree---instead of delays---and 
seeks
to infer the routing tree and packet drop probabilities on the links. This problem
was formalized in~\cite{CaDuHoTo:99} where a maximum-likelihood estimation procedure
was analyzed. In~\cite{CaDuHoTo:99}, the network topology is assumed known.
In more recent independent work, Liang et al.~\cite{LiMoYu:u} (unpublished) 
applied phylogenetic techniques to the inference of the routing topology in this context.
Indeed, the multicast link loss problem is in essence a special case of the standard
model of DNA evolution used in biology.
Similarly to the current paper, Liang et al.~use distance-based techniques. More precisely, 
they give
a computationally efficient reconstruction algorithm with sample complexity $O(b^{-2} \log n)$ where
$b$ (possibly depending on $n$) is a lower bound on the link loss probability. 
Ni and Tatikonda~\cite{NiTatikonda:06,NiTatikonda:07,NiTatikonda:08,NXTY:08} (see above) also considered 
the link loss inference problem.

\subsection{Organization of the Paper}

The paper is organized as follows.
We start with some phylogenetic background in Section~\ref{sec:phylo}.
Our results concerning the topology reconstruction
can be found in the Section~\ref{sec:top}. 
We then present and analyze our delay inference algorithm in
Section~\ref{sec:delay}. 


\section{Phylogenetic Reconstruction Techniques}\label{sec:phylo}

In this section, we summarize and adapt to our setting the DMR algorithm of~\cite{DaMoRo:09}.

\subsection{Basics}

We begin with a few basic notions from phylogenetics.

\paragraph{Tree metrics.} In phylogenetics, the notion of a tree metric is useful 
for reconstructing the topology of phylogenies. We use the notation
$\real_{++} = \{x\in \real:x>0\}$.
\begin{definition}[Tree Metric]\label{def:treemetric}
Let $L$ be a finite set with cardinality $n$. A function $W: L\times L \to \real_{+}$
defines a (nondegenerate) tree metric if the following holds.  There exist a tree
$T=(V,E)$ with leaf set $L$ and a weight function $w: E \to \real_{++}$
such that $W(a,b) = \sum_{e \in P_{ab}} w_e$ for all $a,b \in L$
where $P_{ab}$ is the path between $a$ and $b$ in $T$.
\end{definition}
Tree metrics are usually estimated from samples of the tree process
at the leaves. In that context, Azuma's inequality is useful (see, e.g., \cite{MotwaniRaghavan:95}).
\begin{lemma}[Azuma-Hoeffding Inequality]\label{lemma:azuma}
Suppose $X=(X_1,\ldots,X_k)$ are independent random variables taking values in a set
$S$, and $f:S^k \to \real$ is any $t$-Lipschitz function:
$|f(\bfx) - f(\bfy)|\leq t$ whenever $\bfx$ and $\bfy$ differ at just one coordinate. Then,
$\forall \lambda > 0$,
\begin{equation*}
\prob\left[f(X) - \expec[f(X)] \geq \lambda \right]
\leq \exp\left(-\frac{\lambda^2}{2 t^2 k}\right),
\end{equation*}
and
\begin{equation*}
\prob\left[f(X) - \expec[f(X)] \leq - \lambda \right]
\leq \exp\left(-\frac{\lambda^2}{2 t^2 k}\right).
\end{equation*}
\end{lemma}



\paragraph{Bipartitions.}
A useful combinatorial description of a tree $T=(V,E)$ is obtained
by noticing that each edge $e\in E$ of the tree naturally corresponds to
a partition of the leaves $L$ into two subsets (that is, the leaves on
either ``side'' of $e$). Such partitions are called \emph{bipartitions}
and they characterize the tree: it is easy to generate all bipartitions
corresponding to a given tree, and on the other hand, there is a simple efficient iterative
procedure to recover a tree from the set of all of its bipartitions. 
See~\cite{Felsenstein:04, SempleSteel:03} for details.

\subsection{Distorted Metric Algorithms}\label{sec:distorted}

Classical distance-based reconstruction algorithms (that is, those
methods based on tree metrics) such as UPGMA~\cite{SneathSokal:73} 
or Neighbor-Joining (NJ)~\cite{SaitouNei:87}, 
typically make use of \emph{all} pairwise distances
between leaves.
This leads to difficulties because ``long'' distances are more ``noisy''
and require a large number of samples to be accurately estimated.
For instance, in the phylogenetic context, the widely used NJ algorithm is computationally
efficient, but it is known to require \emph{exponentially} many samples---even for
simple linear trees~\cite{LaceyChang:06}.

An important breakthrough was made in~\cite{ErStSzWa:99a} where it was shown
that it was in fact enough to use ``short'' distances to fully recover the tree
under reasonable assumptions. To help understand this result, we need a
notion of tree ``depth.''
Given an edge $e \in E$, the \emph{chord depth} of $e$ is the length (in graph distance) 
of the shortest path between two leaves on which $e$ lies\footnote{Note that unlike~\cite{DaMoRo:09} we use the graph distance in the definition
of chord depth. Because of our assumptions (see below) the two graph and weighted distances are the same
up to a constant factor. Note also that we are using a different definition than~\cite{ErStSzWa:99a}. But again the difference
is only a constant factor.}. That is,
\begin{equation*}
\pathdepth(e) = \min\left\{\dist(u,v)\ :\ u,v \in L, e\in P_{uv}\right\},
\end{equation*}
where $\dist$ is the graph distance on $T$.
We define the \emph{chord depth of a tree $T$} to be the maximum
chord depth in $T$
\begin{equation*}
\pathdepth(T) = \max\left\{\pathdepth(e)\ :\ e\in E\right\}.
\end{equation*}
It is easy to show that $\depth(T) \leq \log_2 n$ if the degree of all internal nodes
is at least 3 (argue by contradiction). In a nutshell, the key insight behind
the results in~\cite{ErStSzWa:99a} is that the diameter and the depth of a tree
behave very differently: even though the diameter can be as large as $O(n)$,
the depth is always $O(\log n)$, in other words, each edge lies on a ``short'' path
between two leaves. Using clever combinatorial arguments, Erd\"os et al.~\cite{ErStSzWa:99a}
showed that one can reconstruct trees with much fewer samples by ignoring those
distances corresponding to paths longer than $O(\log n)$.

More recently, Daskalakis et al.~\cite{DaMoRo:09} relaxed some of the assumptions
in~\cite{ErStSzWa:99a}. In particular, they gave a reconstruction algorithm
based on short distances allowing internal degrees bigger than 3---which is particularly
relevant in the networking context. Their algorithm, which we will call the DMR algorithm,
reconstructs all bipartitions using only distances smaller than a threshold of order
$O(\log n)$. To check that the algorithm works, one only needs to show that such
distances are accurately estimated for a given number of samples.
In the tomography setting, the DMR algorithm will allow us to reconstruct 
the routing tree using as few as $\poly\log n$ samples (see next section).
The details of the algorithm are sketched in Appendix~\ref{section:dmr}.

We now state a corollary of~\cite{DaMoRo:09} that will be useful to us.
We first need the following definition which formalizes the idea
that short distances are accurately estimated (and that long distances can in some sense be ignored).
\begin{definition}[Distorted Metric~\cite{Mossel:07,KiZhZh:03}]\label{def:distorted metric}
Let $T = (V,E)$ be a tree
with leaf set $L$ and edge weight function $w: E \to \real_{++}$.
Let $W: L\times L \to \real_{+}$ be the corresponding tree metric.
Fix $\tilde\tau, \widetilde M > 0$. 
We say that $\widehat{W} : L\times L \to (0,+\infty]$ is a $(\tilde\tau, \widetilde M)$-\emph{distorted metric}
for $T$ or a $(\tilde\tau, \widetilde M)$-\emph{distortion} of $W$ if:
\begin{enumerate}
\item (Symmetry)
For all $u,v \in L$, $\widehat{W}$ is symmetric, that is,
\begin{equation*}
\widehat{W}(u,v) = \widehat{W}(v,u);
\end{equation*}
\item (Distortion) $\widehat{W}$ is accurate on
``short'' distances, that is, for all $u,v \in L$, if either 
$W(u,v) < \widetilde M + \tilde\tau$ or $\widehat{W}(u,v) < \widetilde M + \tilde\tau$ then
\begin{equation*}
\left|W(u,v) - \widehat{W}(u,v)\right| < \tilde\tau.
\end{equation*}
\end{enumerate}
\end{definition}
Let $f,g > 0$ be bounds
on the edge weights, that is, $f \leq w_e \leq g$ for all $e\in E$.
We say that such an edge weight function satisfies the $(f,g)$-condition.
\begin{theorem}[DMR Algorithm~\cite{DaMoRo:09}]\label{thm:dmr}
Let $ 0 < f < g < +\infty$,
$\tilde\alpha < 1/6$,
and
$\tilde\beta > 2$.
There is a polynomial-time algorithm $\mathbb{A}$ 
such that, for all trees $T = (V,E)$ with edge weight function $w$ satisfying the $(f,g)$-condition
and all $(\tilde\alpha f, \tilde\beta g \depth(T))$-distortions $\widehat W$ of $W$ (where $W$ is the tree metric corresponding to $w$),
$\mathbb{A}$ applied to $\widehat W$ returns $T$.
\end{theorem}
Note that the previous theorem is a deterministic statement about distorted metrics. We show how to
estimate such a distorted metric from random samples with high probability in Section~\ref{sec:samples}.

\section{Routing Tree Reconstruction}\label{sec:top} 

The goal of this section is to
reconstruct efficiently the topology of the routing tree
using Theorem~\ref{thm:dmr}.

\subsection{Variance Metric}

From Definition~\ref{def:treemetric}, one can define a tree metric by
first choosing a tree---in our case, the routing tree---and
then defining a weight function on its edges.
Any positive quantity can serve as a weight. The important point
is that one must be able to estimate the resulting tree metric from samples
at the leaves. This governs the choice of the weight function.

Let $T=(V,E)$ be the (unknown) routing tree with leaf set $L$
and consider the choice of
weights
\begin{equation*}
w^{(2)}_e = \var[d_e],
\end{equation*}
for all $e\in E$ and the corresponding tree metric
\begin{equation*}
W^{(2)}(a,b) \equiv \sum_{e\in P_{ab}} \var[d_e],
\end{equation*}
for all $a,b \in L$.
Our first task is to check that this metric can be estimated from
samples at the leaves. Let $a,b$ be leaves and consider the
quantity $\delta^{(2)}_{ab}\equiv \var[D_a - D_b]$
(where recall from (\ref{eq:total}) that $D_u$ is the delay at $u$). The delays $D_a$ and $D_b$ are observed
at the leaves $a$ and $b$ respectively and therefore the variance of
$D_a - D_b$ can be easily estimated.
Moreover, we claim that the equality $\delta_{ab}^{(2)} = W^{(2)}(a,b)$ holds.
Indeed, denote $\gamma_{ab}$ the common ancestor of $a$ and $b$, that is,
the node at which all three paths $P_{ab}$, $P_{0 a}$, and $P_{0 b}$
intersect (where we assume $a,b\neq 0$). Then, by independence of the edge delays, we have
\begin{equation*}
\delta^{(2)}_{ab}
=\var[D_a - D_b]
= \var\left[\sum_{e\in P_{a \gamma_{ab}}} d_e - \sum_{e\in P_{\gamma_{ab} b}} d_e\right]
= \sum_{e\in P_{a \gamma_{ab}}} \var[d_e] + \sum_{e\in P_{\gamma_{ab}} b} \var[d_e]
= W^{(2)}(a,b).
\end{equation*}
Therefore, we can estimate $W^{(2)}$ by estimating $\delta^{(2)}$ at the leaves.

To estimate $\delta_{ab}^{(2)}$ from $k$ samples, we use the standard unbiased estimator for the variance
of $D_a - D_b$
\begin{equation*}
\hat\delta_{ab}^{(2)} = 
\frac{1}{k-1} \sum_{i = 1}^k \left[(D_a^i - D_b^i) - {\hat \delta}_{ab}^{(1)}\right]^2,
\end{equation*}
where
\begin{equation*}
{\hat \delta}_{ab}^{(1)} = \frac{1}{k}\sum_{i=1}^k (D_a^i - D_b^i).
\end{equation*}
Below, we will need to show that $\hat\delta_{ab}^{(2)}$ is well concentrated around
$\delta_{ab}^{(2)}$, which follows from the Azuma-Hoeffding inequality (see Lemma~\ref{lemma:azuma}). The next lemmas provide the necessary Lipschitz condition.
\begin{lemma}\label{lemma:variance}
Suppose $X=(X_1,\ldots,X_k)$ are independent random variables taking values in
$[-B,B]$ with $k \geq 2$.
Then, the variance estimator
\begin{equation*}
s_X^2 = \frac{1}{k-1} \sum_{i = 1}^k (X_i - \overline X)^2 = \frac{1}{k(k-1)} \sum_{i<j} (X_i - X_j)^2,
\end{equation*}
where $\overline X$ is the sample average, is $\frac{4 B^2}{k}$-Lipschitz.
\end{lemma}
\begin{proof}
Let $X$ be as above and let $Y$ differ from $X$ in one coordinate. Then
\begin{eqnarray*}
\left|s_X^2 - s_Y^2\right|
&\leq&  \frac{1}{k(k-1)} \sum_{i<j} \left|(X_i - X_j)^2 - (Y_i - Y_j)^2\right|
\leq \frac{1}{k} 4 B^2.
\end{eqnarray*}
\end{proof}
We then get immediately the following.
\begin{lemma}[Lipschitz Constant for Delay-based Metric]\label{lemma:lip}
Say $\hat\delta_{ab}^{(2)}$ is computed with $k$ samples.
Then, $\hat\delta_{ab}^{(2)}$
is then $\frac{4 |P_{ab}|^2 M^2}{k}$-Lipschitz.
\end{lemma}
$\blacksquare$

\subsection{Inferring the routing tree}\label{sec:samples}

Equipped with a legitimate tree metric, we use the DMR algorithm to infer the
topology. Here, we use Theorem~\ref{thm:dmr} to prove that
the routing tree can be inferred with $\poly\log n$ samples at the leaves.
This is our main result for this section. The main technical
difficulty (unlike the phylogenetic case) is in controlling the deviation of
``long distances.'' (See second part of the proof.)
Fix $\tilde \alpha$, $\tilde \beta$, $f$, $g$ as in Theorem~\ref{thm:dmr}.
Note that by assumption we have
\begin{equation}\label{eq:assumpw2}
f \leq w^{(2)}_e \leq g,
\end{equation}
for all $e$ with
\begin{equation*}
g \equiv M^2.
\end{equation*}
\begin{theorem}[Efficient Network Inference]\label{thm:rti} 
Let $T=(V,E)$ be the (unknown)
routing tree where edge delays satisfy
Assumptions~\ref{assump:regularity} and~\ref{assump:moments}. 
Consider the tree metric $W^{(2)} =
\delta^{(2)}$ and assume that the estimate ${\widehat{W}}^{(2)} =
{\hat\delta}^{(2)}$ is computed using $k$ samples at the leaves.
Then, DMR returns the correct topology for $T$ with probability
$1- o(1)$ if $k = \Omega(\log^5 n)$ 
(where the constant factor depends only on $f,g$), as $n$ tends to $+\infty$.
\end{theorem}
\begin{proof}
Assume $k$ is as stated above.
We apply Theorem~\ref{thm:dmr} and therefore only need to show that 
${\widehat{W}}^{(2)}$ is a $(\tilde\alpha f, \tilde\beta g \depth(T))$-distortion of
$W^{(2)}$ when $k = \Omega(\log^5 n)$.

\paragraph{Part 1.} First, we must show that distances smaller than $\tilde\beta g \depth(T) + \tilde\alpha f$
under $W^{(2)}$ are approximated within $\tilde\alpha f$. For reasons that will become
clear below, we show instead that distances smaller than \emph{twice} that amount are well approximated.
Let $a,b$ be
a pair of leaves at distance at most $2\tilde\beta g \depth(T) + 2\tilde\alpha f$. 
Let $\acal$ be the probability that, for all such pairs, $W^{(2)}$ is approximated within
$\tilde\alpha f$.
By our assumption (\ref{eq:assumpw2}), the number of edges
on the path between $a$ and $b$ is at most
\begin{equation*}
|P_{ab}| \leq (2\tilde\beta\depth(T) + 1)\frac{g}{f},
\end{equation*} 
where we used that $f < g$ and $\tilde\alpha < 1/6$.
By
Lemmas~\ref{lemma:azuma} and~\ref{lemma:lip}, we have
\begin{eqnarray*}
\prob\left[\left|\delta_{ab}^{(2)} - \hat\delta_{ab}^{(2)}\right| \geq \tilde\alpha f\right]
&\leq& 2\exp\left(-\frac{(\tilde\alpha f)^2 k}{2 [4 (2\tilde\beta\depth(T) + 1)^2\frac{g^2}{f^2} M^2]^2 }\right)
\leq \frac{1}{\poly(n)},
\end{eqnarray*}
from $\depth(T) = O(\log n)$, $k = \Omega(\log^5 n)$, and the fact that $f,g,M$ are constants.
The notation $\poly(n)$ means $O(n^K)$ for a $K$ {\it as a large as we need} as
long as the constant factor in $k$ is large enough.
Since there are at most $n^2$ such pairs of leaves, we get $\acal \leq \frac{1}{\poly(n)}$.

\paragraph{Part 2.} Let $a,b$ be
a pair of leaves at distance at least $2\tilde\beta g \depth(T) + 2\tilde\alpha f$
under $W^{(2)}$. 
We now show that, for all such pairs, ${\widehat{W}}^{(2)}$ is at least $\tilde\beta g \depth(T) + \tilde\alpha f$.
Let $\bcal$ be the probability of that event. 
Note first that from Azuma-Hoeffding (Lemma~\ref{lemma:azuma}),
it follows that for any pair of leaves $a,b$,
\begin{equation}\label{eq:sqrt}
\prob\left[\left|
(D_a - D_b) - \expec[D_a - D_b]
\right| \leq \sqrt{|P_{ab}|}\ \Theta(\sqrt{ \log n})\right] \geq 1 - \frac{1}{\poly(n)}.
\end{equation}
Let $\ecal$ be the event that the inequality in square brackets in (\ref{eq:sqrt}) holds for all
$k$ samples used to compute $\hat\delta_{ab}^{(2)}$. Then from
Lemma~\ref{lemma:variance}, on $\ecal$,
the Lipschitz constant of $\hat\delta_{ab}^{(2)}$ (as a function of the centered
samples $(D^i_a - D^i_b) - \expec[D_a - D_b]$) is 
$t = \frac{|P_{ab}|}{k}\Theta\left(\log n\right)$ and therefore, by Lemma~\ref{lemma:azuma} again,
\begin{eqnarray*}
\prob\left[\hat\delta_{ab}^{(2)} \leq \tilde\beta g \depth(T) + \tilde\alpha f
\ \bigg|\ \ecal\ \right]
&\leq& \prob\left[\hat\delta_{ab}^{(2)} \leq \frac{\expec[\hat\delta_{ab}^{(2)}|\ecal]}{2}\ \bigg|\ \ecal\ \right]\\
&\leq& \prob\left[\expec[\hat\delta_{ab}^{(2)}|\ecal] - \hat\delta_{ab}^{(2)} \geq \frac{\expec[\hat\delta_{ab}^{(2)}|\ecal]}{2}
\ \bigg|\ \ecal\ \right]\\
&\leq& \exp\left(-\frac{(\expec[\hat\delta_{ab}^{(2)}|\ecal]/2)^2}{2 t^2 k}\right)\\
&\leq& \exp\left(- \frac{k}{O(\log^2 n)}\right)\\
&\leq& \frac{1}{\poly(n)},
\end{eqnarray*}
where we used $\delta_{ab}^{(2)} =  \Theta(|P_{ab}|)$ and
\begin{equation*}
(1-o(1))\delta_{ab}^{(2)} \leq \expec[\hat\delta_{ab}^{(2)}|\ecal] \leq (1+o(1))\delta_{ab}^{(2)},
\end{equation*}
which follows from $\expec[\hat\delta_{ab}^{(2)}|\ecal]\prob[\ecal]
+ \expec[\hat\delta_{ab}^{(2)}|\ecal^c]\prob[\ecal^c] = \expec[\hat\delta_{ab}^{(2)}]$, 
$\expec[\hat\delta_{ab}^{(2)}] = \delta_{ab}^{(2)}$, $\prob[\ecal^c] \leq 1/\poly(n)$, and
$\hat\delta_{ab}^{(2)} = O(n^2)$. Therefore, we have $\bcal \leq 1/\poly(n)$.

Combining the two parts of the argument, we have shown that, except with $o(1)$ probability,
${\widehat{W}}^{(2)}$ is a $(\tilde\alpha f, \tilde\beta g \depth(T))$-distortion of
$W^{(2)}$. Indeed, by Part 2 the pairs of leaves for which ${\widehat{W}}^{(2)} < \tilde\beta g \depth(T) + \tilde\alpha f$
must have $W^{(2)} < 2\tilde\beta g \depth(T) + 2\tilde\alpha f$ and such pairs satisfy the approximation
guarantee required by the definition of a distorted metric by Part 1. Moreover, Part 1 
implies in particular that pairs of leaves such that $W^{(2)} < \tilde\beta g \depth(T) + \tilde\alpha f$
also satisfy the approximation guarantee.
\end{proof}

\section{Edge Delay Inference}\label{sec:delay}

In this section, we show how to
estimate the characteristic moments of edge delays. 
In Section~\ref{sec:top}, we showed how to reconstruct the 
topology efficiently with high probability (see Theorem~\ref{thm:rti}).
Therefore, along with Assumptions~\ref{assump:regularity}
and~\ref{assump:moments}, we make the following assumption.
\begin{assumption}[Correct Reconstruction of Routing Tree]
\label{assump:tree correctly est} We assume that the routing tree
was {\em correctly} estimated.~(This is true with high probability
by Theorem~\ref{thm:rti}.)
\end{assumption}

Our general idea to recover delay distributions 
is to define so-called ``additive functions'' whose edge weights
are moments of delays. Then we use the AFI algorithm below
to recover the moments efficiently from the data at the leaves. 
As it turns out, even moments are rather straightforward to estimate inductively while
odd moments are trickier.
Also, as in the tree reconstruction algorithm (see also~\cite{ErStSzWa:99a,MosselRoch:06}),
the AFI algorithm uses only ``short''  paths
during the estimation process, which allows a significant reduction in the
sample size (see Propositions~\ref{prop:concineq},~\ref{prop:concineqgen}
 and Theorems~\ref{thm:symer},~\ref{thm:generaler} 
for details).

\subsection{Additive Functions}

In the remainder of this paper, we use additive metric-type ideas to
estimate moments of edge delays. For this purpose,
we need to recover edge weights from appropriately defined
tree metrics.
In fact, we use a notion of ``generalized'' tree metric which
is useful in treating odd moments. This definition allows
for negative edge weights.
\begin{definition}[Additive function]\label{def:additive function} 
A function on the leaf set of the
tree $W:L\times L \to \real$ is called an {\em additive function} on the
leaves if there exists weights $w_e\in \real$ on each of the edges (not
necessarily positive), such that for all leaves $a,b$
\begin{equation*}
W(a,b)= \sum_{e \in P_{ab}} w_e.
\end{equation*}
\end{definition}

\begin{figure}
\begin{center}
\input{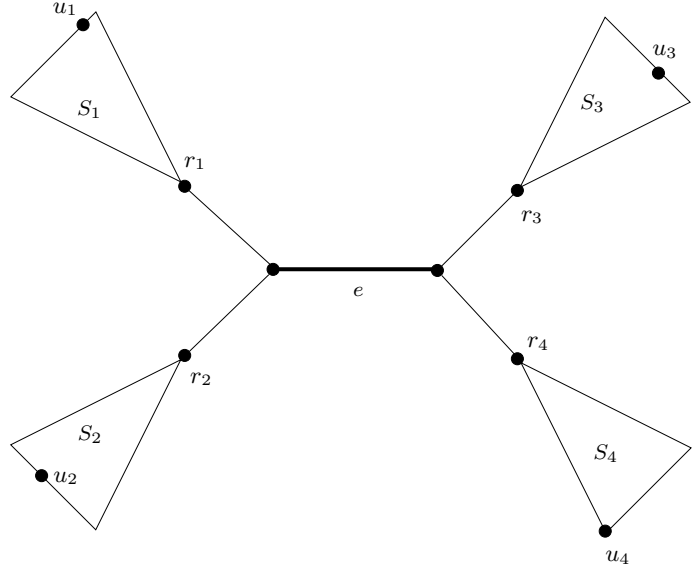}
\caption{Edge weight inference.}\label{fig:quartet4} 
\end{center}
\end{figure}

Suppose we are given access to an additive function
$W$ on the leaves. Our goal is now 
to recover the $w_e$'s from the function $W$, assuming further
that we are given the tree $T$.
For this purpose, we use a standard algorithm from combinatorial phylogenetics---related
to the so-called Four-Point Method of Buneman~\cite{Buneman:71} (see also~\cite{Felsenstein:04,SempleSteel:03}).
We will refer to this algorithm as the 
\textsc{Additive Function Inference} (AFI) algorithm.
See Figures~\ref{fig:afi} and~\ref{fig:quartet4}.
\begin{figure}[!ht]
\framebox{
\begin{minipage}{16cm}
{\small
\textbf{Algorithm} \textsc{Additive Function Inference} \\
\textit{Input:} tree $T$, function $W$ at the leaves;\\
\textit{Output:} edge weights $w_e$, for all $e\in E$;

\begin{itemize}
\item For all internal edges $e$,
\begin{itemize}
\item Let $S_1,\ldots,S_4$ be the four subtrees hanging from
$e$ as in Figure~\ref{fig:quartet4};
\item For each $S_i$, compute $u_i$  the closest (in graph
distance) leaf
to the root $r_i$ of $S_i$;
\item Compute
\begin{equation*}
w_e=\frac{1}{2}(W(u_1,u_3)+W(u_2,u_4)-W(u_1,u_2)-W(u_3,u_4)).
\end{equation*}
\end{itemize} 
\item For all leaf edges $e$,
\begin{itemize}
\item Let $e = (a,v)$ with $a$ a {\em leaf}; 
\item Proceed as above where $u_3$ and $u_4$ are set to $a$.
\end{itemize}
\end{itemize}

}
\end{minipage}
}
\caption{Algorithm \textsc{Additive Function Inference}.} \label{fig:afi}
\end{figure}

\subsection{Delay-based metrics}

Let $T=(V,E)$ be the routing tree with leaf set $L$
and consider again the choice of
weights
\begin{equation*}
w^{(2)}_e = \var[d_e],
\end{equation*}
for all $e\in E$ and
\begin{equation*}
W^{(2)}(a,b) = \sum_{e\in P_{ab}} w^{(2)}_e,
\end{equation*}
for all $a,b \in L$.
Recall that
\begin{equation*}
\delta^{(2)}_{ab}
=\var[D_a - D_b]
= \var\left[\sum_{e\in P_{a \gamma_{ab}}} d_e - \sum_{e\in P_{\gamma_{ab} b}} d_e\right]
= \sum_{e\in P_{a \gamma_{ab}}} \var[d_e] + \sum_{e\in P_{\gamma_{ab}} b} \var[d_e]
= W^{(2)}(a,b).
\end{equation*}
Therefore, using the AFI algorithm, we can recover estimates of the $w^{(2)}_e$'s from 
the $\hat\delta^{(2)}_{ab}$'s. 

More generally, we let 
\begin{equation*}
w^{(j)}_e = \expec\left[\left({\bar d}_e\right)^j\right],
\end{equation*}
for all $e\in E$ where
\begin{equation*}
{\bar d_e} = d_e - \expec\left[d_e\right].
\end{equation*} 
Also, let
\begin{equation*}
W^{(j)}(a, b) = \sum_{e\in P_{ab}} w^{(j)}_{e},
\end{equation*}
for all $a,b \in L$. Let
\begin{equation*}
\overline{D}_a = D_a - \expec\left[D_a\right],
\end{equation*}
for all $a \in L$.
Again, to obtain $W^{(j)}(a,b)$, we seek to use the quantity
\begin{equation*}
\delta^{(j)}_{ab} = \expec\left[\left(\overline{D}_a - \overline{D}_b\right)^j\right],
\end{equation*}
for $j > 1$, which can be estimated from the samples using
\begin{equation*}
\hat{\delta}_{ab}^{(j)}
=\frac{1}{k} \sum_{i=1}^{k} \left(\left(D_a^i-D_b^i\right) - {\hat\delta}_{ab}^{(1)}\right)^j
 \end{equation*}
 where
 \begin{equation*}
 {\hat \delta}_{ab}^{(1)} = \frac{1}{k}\sum_{i=1}^{k} \left(D_a^i - D_b^i\right).
 \end{equation*}
As Lemma~\ref{lem:even} below shows, this can be done inductively. However, the lemma
also shows that odd moments have to be treated more carefully.

\subsection{Algorithm for Moment Inference}\label{sec:er}

We first need the following
definitions. Let $a,b$ be leaves and $j\in \nintgr$. 
We use the notation $[h] = \{0,\ldots,h\}$ for $h \in\nintgr$. 
Recall that $\gamma_{ab}$ is the most recent common ancestor
of $a$ and $b$ in the tree. Denote $\nu = |P_{ab}|$, $\alpha = |P_{a \gamma_{ab}}|$,
and $\beta = |P_{\gamma_{ab}b}|$, and
define
\begin{equation*}
\dcal_j(a,b) 
= \left\{
(\bfx, \bfy)\in[j-1]^{\alpha}\times[j-1]^{\beta}\,:\,
\sum_{i=1}^\alpha x_i + \sum_{i=1}^\beta y_i = j
\right\}.
\end{equation*}
For $(\bfx,\bfy) \in \dcal_j(a,b)$, let
\begin{equation*}
\binom{j}{\bfx,\bfy} 
=
\frac{j!}{\prod_{i=1}^\alpha x_i! \prod_{i=1}^\beta y_i!}.
\end{equation*}
and consider the function
\begin{equation*}
\fcal_j(a,b) 
=
\sum_{(\bfx,\bfy)\in\dcal_j(a,b)} \binom{j}{\bfx,\bfy} \prod_{i=1}^\alpha w_{e_i}^{(x_i)}
\prod_{i=1}^\beta (-1)^{y_i }w_{f_i}^{(y_i)},
\end{equation*}
where $P_{a\gamma_{ab}} = (e_1,\ldots,e_\alpha)$
and
$P_{\gamma_{ab}b} = (f_1,\ldots,f_\beta)$.
\begin{lemma}\label{lem:even}
Let $j\in\nintgr$ and define the function $\fcal_{j} :L \times L \rightarrow \real $ as above.
Then, 
\begin{enumerate}
\item we have for all $a,b \in L$
\begin{equation}\label{eq:wj} 
\delta_{ab}^{(j)}
- \fcal_{j}(a,b)
= 
\sum_{i=1}^\alpha w_{e_i}^{(j)} + (-1)^j \sum_{i=1}^\beta w_{f_i}^{(j)},
\end{equation}
\item in particular, if $j$ is even, we have for all $a,b \in L$
\begin{equation}\label{eq:weven}
\delta_{ab}^{(j)}
- \fcal_{j}(a,b) = W^{(j)}(a,b).
\end{equation}
\end{enumerate}
\end{lemma}
\begin{proof}
This follows immediately from a multinomial expansion.
\end{proof}
The important point to note in (\ref{eq:weven}) is that $\fcal_j(a,b)$ depends only 
on delay moments {\em of order strictly less than $j$} and that 
$\delta_{ab}^{(j)}$ can be estimated from samples at the leaves. 
Therefore, if $j$ is even and if we have estimates of all edge delay moments
of order up to $j-1$, we can estimate $W^{(j)}(a,b)$ by~(\ref{eq:weven}). Using the AFI algorithm,
we can then get an estimate of the $j$-th moments $w_e^{(j)}$. However, if $j$ is odd, the coefficient
$(-1)^j$ in~(\ref{eq:wj}) precludes the use of this procedure. Lemma~\ref{lem:odd} below shows
how to handle this case. We note in passing that Lemma~\ref{lem:even} above is sufficient
for delay distributions symmetric about their mean. 
Indeed, in that case, all odd central moments are zero and one can use~(\ref{eq:weven}) recursively
to estimate all even characteristic moments.
See Figure~\ref{fig:symer}. 
\begin{figure}[!ht]
\framebox{
\begin{minipage}{16cm}
{\small
\textbf{Algorithm} \textsc{Symmetric Edge Reconstruction} \\
\textit{Input:} data $\{D^1_a\}_{a\in L}, \ldots, \{D^k_a\}_{a\in L}$ at the leaves; topology $T$;\\
\textit{Output:} estimated characteristic (even) moments ${\hat w}_e^{(j)}$ for all $e\in E$
and $2\leq j\leq J$ even;

\begin{itemize}
\item Initialization: set all estimates of odd moments to 0;
\item Main Loop: For all $2\leq j\leq J$ even,
\begin{itemize}
\item For all $a,b \in L$,
\begin{itemize}
\item Estimate ${\hat \delta}_{ab}^{(j)}$;
\item Estimate $\fcal_j(a,b)$ with
\begin{equation*}
\widehat{\fcal}_j(a,b) 
=
\sum_{(\bfx,\bfy)\in\dcal_j(a,b)} \binom{j}{\bfx,\bfy} \prod_{i=1}^\alpha {\hat w}_{e_i}^{(x_i)}
\prod_{i=1}^\beta (-1)^{y_i }{\hat w}_{f_i}^{(y_i)}.
\end{equation*} 
\item Compute
\begin{equation*}
\widehat{W}^{(j)}(a,b)
=
{\hat \delta}_{ab}^{(j)}
- \widehat{\fcal}_{j}(a,b)
\end{equation*}
\end{itemize}
\item Use the AFI algorithm on $\widehat{W}^{(j)}(a,b)$ to recover
all ${\hat w}_e^{(j)}$'s.
\end{itemize}
\end{itemize}
}
\end{minipage}
}
\caption{Algorithm \textsc{Symmetric Edge Reconstruction}.} \label{fig:symer}
\end{figure}

We now tackle odd moments. A proper estimation procedure follows from the next lemma.
We first need a few definitions. For $a,b \in L$, and $1\leq i^*\leq \alpha$, we let
\begin{equation*}
\ecal^{(1)}_j(a,b;i^*) 
= \left\{
(\bfx, \bfy)\in[j-1]^{\alpha}\times[j-1]^{\beta}\,:\,
\sum_{i=1}^\alpha x_i + \sum_{i=1}^\beta y_i = j,\ x_{i^*} \geq 1
\right\}.
\end{equation*}
and
\begin{equation*}
\gcal^{(1)}_j(a,b) 
=
\sum_{i^* = 1}^\alpha\sum_{(\bfx,\bfy)\in\ecal^{(1)}_j(a,b;i^*)}  x_{i^*} \binom{j-1}{\bfx,\bfy} 
\prod_{i=1}^\alpha w_{e_i}^{(x_i)}
\prod_{i=1}^\beta (-1)^{y_i }w_{f_i}^{(y_i)},
\end{equation*}
where we use the notations of Lemma~\ref{lem:even}.
Similarly, for $1\leq i^*\leq \beta$, we define 
$\ecal_j^{(2)}(a,b;i^*)$ and $\gcal_j^{(2)}(a,b)$ by interchanging the roles of $\bfx$
and $\bfy$.
Our next definition requires a few combinatorial notions. 
Recall the definition of quartet split from Section~\ref{sec:phylo}.
Let 
$a,b,c$ be any leaves in a {\it rooted} tree $T$ with root $0$
(which is also a leaf). 
We write $ab|c$
if $ab|c0$ holds in $T$.
Then, for all leaves $a,b,c\neq 0$ with $ab|c$, let
\begin{equation*}
\phi_{ab|c}^{(j)}
=
\expec\left[\left(\overline{D}_a - \overline{D}_b\right)^{j-1}
\left(\overline{D}_a + \overline{D}_b - 2\overline{D}_c\right)
\right].
\end{equation*}
\begin{lemma}\label{lem:odd}
Let $j\in\nintgr$.
Then, using the notations above,
we have for all $a,b,c \in L$
\begin{equation}\label{eq:wodd}
W^{(j)}(a,b)
=
\phi_{ab|c}^{(j)}
- \left[\gcal^{(1)}_{j}(a,b) + \gcal^{(2)}_{j}(a,b)\right].
\end{equation}
\end{lemma}
\begin{proof}
We write
\begin{equation*}
\expec\left[\left(\overline{D}_a - \overline{D}_b\right)^{j-1}
\left(\overline{D}_a + \overline{D}_b - 2\overline{D}_c\right)
\right]
= \expec\left[\left(\overline{D}_a - \overline{D}_b\right)^{j-1}
\left(\overline{D}_a - \overline{D}_c\right)
\right]
+
\expec\left[\left(\overline{D}_a - \overline{D}_b\right)^{j-1}
\left(\overline{D}_b - \overline{D}_c\right)
\right]
\end{equation*}
Let (as in Figure~\ref{fig:oddrecursion})
\begin{eqnarray*}
H_1 
= \sum_{e\in P_{a \gamma_{ab}}} {\bar d}_e\qquad
H_2 
= \sum_{e\in P_{b \gamma_{ab}}} {\bar d}_e\qquad
H_3 
= \sum_{e\in P_{\gamma_{ac} \gamma_{ab}}} {\bar d}_e\qquad 
H_4 
= \sum_{e\in P_{c \gamma_{ac}}} {\bar d}_e.
\end{eqnarray*}
Note that all these random variables are independent and have
$0$ mean.
Then
\begin{eqnarray*}
\expec\left[\left(\overline{D}_a - \overline{D}_b\right)^{j-1}
\left(\overline{D}_a - \overline{D}_c\right)
\right]
&=&
\expec\left[\left(H_1 - H_2\right)^{j-1}
\left(H_1 + H_3 - H_4\right)
\right]\\
&=& \expec\left[\left(H_1 - H_2\right)^{j-1}
\left(H_1\right)
\right]\\
&=& 
\expec\left[\left(\sum_{e\in P_{a \gamma_{ab}}} {\bar d}_e -  \sum_{e\in P_{b \gamma_{ab}}} {\bar d}_e\right)^{j-1}
\left(\sum_{e\in P_{a \gamma_{ab}}} {\bar d}_e\right)
\right]\\
&=&
\sum_{e\in P_{a\gamma_{ab}}} w_e^{(j)}
+ \gcal_j^{(1)}(a,b).
\end{eqnarray*}
Similarly,
\begin{eqnarray*}
\expec\left[\left(\overline{D}_a - \overline{D}_b\right)^{j-1}
\left(\overline{D}_b - \overline{D}_c\right)
\right]
&=&
\sum_{e\in P_{b\gamma_{ab}}} w_e^{(j)}
+ \gcal_j^{(2)}(a,b).
\end{eqnarray*}
The result follows.
\begin{figure}
\begin{center}
\input{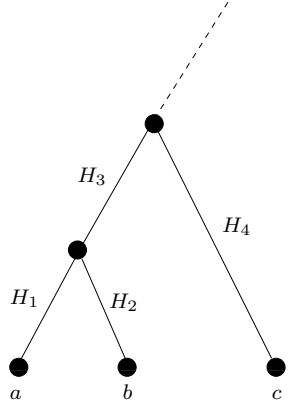}
\caption{The $H_i$'s are centered sums of delays on the corresponding paths.}\label{fig:oddrecursion}
\end{center}
\end{figure}
\end{proof}
Again, the key point in (\ref{eq:wodd}) is that
$\gcal_j^{(1)}(a,b)$ and $\gcal_j^{(2)}(a,b)$ depend only 
on moments {\em of order strictly less than $j$} and that 
$\phi_{ab|c}^{(j)}$ can be estimated from samples at the leaves.
The algorithm for the general case is detailed in Figure~\ref{fig:er}.
We use the plugin estimator for $\phi_{ab|c}^{(j)}$,
\begin{equation*}
{\hat \phi}_{ab|c}^{(j)}
=\frac{1}{k} \sum_{i=1}^{k} \left(\left(D_a^i-D_b^i\right) - {\hat\delta}_{ab}^{(1)}\right)^{j-1}
\left(\left(D_a^i-D_c^i\right) - {\hat\delta}_{ac}^{(1)} 
+ \left(D_b^i-D_c^i\right) - {\hat\delta}_{bc}^{(1)}\right).
\end{equation*}
\begin{figure}[!ht]
\framebox{
\begin{minipage}{16cm}
{\small
\textbf{Algorithm} \textsc{Edge Reconstruction} \\
\textit{Input:} data $\{D^1_a\}_{a\in L}, \ldots, \{D^k_a\}_{a\in L}$ at the leaves; topology $T$;\\
\textit{Output:} estimated characteristic moments ${\hat w}_e^{(j)}$ for all $e\in E$
and $2\leq j\leq J$;

\begin{itemize}
\item Initialization: set all estimates of first moments to 0;
\item Main Loop: For all $2\leq j\leq J$,
\begin{itemize}
\item For all $a,b \in L$,
\begin{itemize}
\item Pick the closest leaf $c$ above $\gamma_{ab}$
\item Compute ${\hat \phi}_{ab|c}^{(j)}$, the plug-in estimator for $\phi_{ab|c}^{(j)}$;
\item Estimate $\gcal^{(1)}_j(a,b)$ with
\begin{equation*}
\widehat{\gcal}^{(1)}_j(a,b) 
=
\sum_{i^* = 1}^\alpha\sum_{(\bfx,\bfy)\in\ecal^{(1)}_j(a,b;i^*)}  x_{i^*} \binom{j-1}{\bfx,\bfy} 
\prod_{i=1}^\alpha {\hat w}_{e_i}^{(x_i)}
\prod_{i=1}^\beta (-1)^{y_i }{\hat w}_{f_i}^{(y_i)},
\end{equation*} 
and similarly for $\gcal^{(2)}_j(a,b)$;
\item Compute
\begin{equation*}
\widehat{W}^{(j)}(a,b)
=
{\hat \phi}_{ab|c}^{(j)}
- \left(\widehat{\gcal}^{(1)}_{j}(a,b)
+ \widehat{\gcal}^{(2)}_{j}(a,b)\right),
\end{equation*}
\end{itemize}
\item Use the AFI algorithm on $\widehat{W}^{(j)}(a,b)$ to recover
all ${\hat w}_e^{(j)}$'s.
\end{itemize}
\end{itemize}
}
\end{minipage}
}
\caption{Algorithm \textsc{Edge Reconstruction}.} \label{fig:er}
\end{figure}

\section{Analysis of the ER Algorithm}\label{section:conc ineq} 

We start with the analysis of the symmetric case.
 
We begin with a concentration result for the estimate ${\hat \delta}_{ab}^{(j)}$.
For convenience, we assume $M\geq 1$. (This can always be obtained
by rescaling.) 
Recall the definition of the depth of $T$ from Section~\ref{sec:phylo} and
remember that $\depth(T) = O(\log n)$ if the degree of all internal nodes
is at least 3.
The dependence of our bounds on the depth of the routing tree
explains the importance of using short paths in the estimation
procedures. 
\begin{proposition}\label{prop:concineq} 
Let $a, b\in L$ 
at graph distance less than $2\depth$ where $\depth = \depth(T)$ is the chord depth
of $T$. Fix $j \in \nintgr$.
We have the following (where the constants depend on $J$ and $M$ only):
\begin{enumerate}
\item
There exists a constant $C$ such that, $\forall \lambda > 0$,
\begin{equation}\label{eqn:delhat} 
\prob\left(\left|\delhat^{(j)}-\expec\left[\delhat^{(j)}\right]\right| >
\lambda\right) \leq 2 \exp\left(-\frac{{\lambda}^2 k}{C \depth^{2j-1}}\right).
\end{equation}

\item There exists a constant $C'$ such that
\begin{equation}\label{eqn:moment bound} 
\expec\left( \left| \delta_{ab}^{(1)} - \delhat^{(1)} \right|^j \right) 
\leq C' \frac{\depth^{j}}{k^{j/2}},
\end{equation}
where $\delta^{(1)}_{ab} = \expec[D_a - D_b]$.

\item There exists a constant $C''$ such that, if $k \geq \depth^2$,
\begin{equation}
\left|\expec\left[\delhat^{(j)}\right] - \del^{(j)}\right| 
\leq C'' \frac{M^{2j}\depth^{j+1}}{\sqrt{k}}.
\end{equation}

\item If further
\begin{equation*}
C'' \frac{M^{2j}\depth^{j+1}}{\sqrt{k}} \leq \lambda,
\end{equation*}
then we have
\begin{equation}
\prob\left[\left|\delhat^{(j)}-\del^{(j)} \right| >
2\lambda\right] \leq 2 \exp\left(-\frac{{\lambda}^2 k}{C \depth^{2j-1}}\right).
\end{equation}
\end{enumerate}
\end{proposition}
\begin{proof}
1. We use Azuma's inequality (see Lemma
\ref{lemma:azuma}). 
Let
\begin{eqnarray*}
\kcal_i &=& (D_a^i - D_b^i)- \delhat^{(1)},
\end{eqnarray*}
where $D_u^i$ is the $i$-th delay sample at node $u$.
Because $|P_{ab}| \leq 2\depth$ and $d_e \in [0,M]$ for 
all $e$, it follows that
\begin{eqnarray*}
|\kcal_i | \leq 4M\depth.
\end{eqnarray*}
Then let
\begin{eqnarray*}\label{eqn:L fn definition}
\lcal &=& \frac{1}{k} \sum_{i=1}^{k} (\kcal_i)^j,
\end{eqnarray*}
and let $\lcal'$ be the same quantity when 
an arbitrary $d^i_e$ is perturbed by $\delta$ with $|\delta| \leq M$
(where $d_e^i$ is the $i$-th delay sample on edge $e$).
Without loss of generality, assume the perturbation
is in the first sample. Then,
\begin{eqnarray}\label{eq:lprime}
\lcal' &=& \frac{1}{k} \left( \left(\kcal_1+ \frac{(k-1)}{k}\delta\right)^j +
\sum_{i=2}^{k} \left(\kcal_i-\frac{\delta}{k}\right)^j \right).
\end{eqnarray}
Now expanding (\ref{eq:lprime}),
we get
\begin{eqnarray*}
|\lcal - \lcal'| 
&\leq&
\frac{1}{k}\left(
2^j(4M\depth)^{j-1} M +
(k-1) \left(2^j (4M\depth)^{j-1} \frac{M}{k}\right)
\right)\\ 
&\leq&
C \frac{\depth^{j-1}}{k},
\end{eqnarray*}
for some constant $C$ depending on $M, J$.
Noting that $\lcal$ depends on at most  $2 \depth k$ random variables $d_e^i$,
we get the result by an application of Azuma's inequality (for a different $C$).

2. Note that
\begin{eqnarray*}
\lcal &=& \del^{(1)} - \delhat^{(1)},
\end{eqnarray*}
is a $\frac{2M \depth}{k}$ -Lipschitz function of $\{D_a^i - D_b^i\}_{i\in[k]}$ 
thus we have by Azuma's inequality
\begin{equation*}
\prob\left[\left| \del^{(1)} - \delhat^{(1)} \right| > \lambda \right]
\leq 2 \exp \left( - \frac{k \lambda^2}{8 M^2 \depth^2} \right).
\end{equation*}
Now we use the fact that for a positive random variable $Y$,
\begin{equation*}
\expec \left[Y^j\right] = j \int_0^{\infty} \lambda^{j-1} \prob(Y > \lambda ) \diff \lambda.
\end{equation*}
If $Y = |\del^{(1)} - \delhat^{(1)}|$ and $\psi =  \frac{k}{8 M^2 \depth^2}$, we have
\begin{eqnarray*}
\expec \left[Y^j\right]
\leq \psi^{-\frac{j}{2}} \int_{0}^{+\infty} y^{\frac{j}{2} - 1} e^{-y}\diff y
= \left( \frac{8 M^2 \depth^2}{k}\right)^{j/2} C'.
\end{eqnarray*}
That proves 2 (for a different $C'$).

3. We have
\begin{eqnarray*}
\delhat^{(j)} &=& \frac{1}{k} \sum_{i=1}^{k} \left(\left(D_a^i -D_b^i \right)- \delhat^{(1)} \right)^j
= \frac{1}{k}\sum_{i=1}^{k} \left(\left(D_a^i -D_b^i- \del^{(1)}\right) + \left(\del^{(1)} - \delhat^{(1)} \right) \right)^j.
\end{eqnarray*}
Now expand using the binomial theorem and take expectations to get
\begin{eqnarray*}
\left|\expec\left[\delhat^{(j)}\right]- \del^{(j)} \right| 
&\leq& \frac{1}{k}\expec \left|
\sum_{i=1}^{k} \sum_{h=0}^{j-1} {j \choose h } \left(D_a^i - D_b^i - \del^{(1)} \right)^h
(\delhat^{(1)} - \del^{(1)} )^{j-h} \right|\\
&\leq& C'' (4M \depth)^j \max_{0\leq h\leq j-1 }\left\{\expec\left|\del^{(1)} - \delhat^{(1)} \right|^{j-h}\right\}.
\end{eqnarray*}
Note that by $k \geq \Delta^2$, it follows that the maximum is attained at $h=j-1$ in (\ref{eqn:moment bound}).

4. This follows from 1. and 3.
\end{proof}

We then get the main theorem in the symmetric case. 
Recall that $J = O(1)$ and that, in general, 
$\depth = O(\log n)$ where $n$ is the number of leaves.
\begin{theorem}\label{thm:symer} 
Let $\eps > 0$ be arbitrarily small. If 
$k = \omega(\depth^{2J^2} \log n)$, then after an application of \textsc{SymER}, one has
\begin{equation}
\prob \left[ \left|{\hat\delta}_e^{(j)}-\delta_e^{(j)}\right| 
\leq \eps,\ 
\forall e\in E,\ \forall 1\leq j \leq J \right] \geq 1-o(1),
\end{equation}
as $n \to +\infty$. The algorithm runs in time $O(\depth^J n^2)$.
\end{theorem}
\begin{proof}
Let $(a,b) \in L\times L$ be called a {\em short pair} if 
$a,b$ are at graph distance at most $2\depth$.
Denote $\scal$ be the set of all short pairs.
Let
\begin{equation*}
\sigma_j 
= 
\max_{(a,b) \in \scal}
\left|
\widehat{W}^{(j)}(a,b) - W^{(j)}(a,b)
\right|,
\end{equation*}
and
\begin{equation*}
\Sigma_j = \max_{1\leq i\leq j} \sigma_i.
\end{equation*}
It follows immediately from the application of the AFI
algorithm that
\begin{equation*}
\max_{e\in E}
\left|
{\hat w}^{(j)}_e - w^{(j)}_e
\right|
\leq 2\sigma_j.
\end{equation*}
Therefore, it suffices to prove
\begin{equation*}
\Sigma_J = o(1),
\end{equation*}
with high probability as $n$ tends to $+\infty$.

Further, assume we have a uniform bound
\begin{equation*}
\max_{1\leq j\leq J} \max_{(a,b) \in \scal}
\left|
\delhat^{(j)} - \del^{(j)}
\right| \leq \tau^*.
\end{equation*}
Recall that
\begin{equation*}
\widehat{W}^{(j)}(a,b)
=
{\hat \delta}_{ab}^{(j)}
- \widehat{\fcal}_{j}(a,b)
\end{equation*}
where
\begin{equation*}
\widehat{\fcal}_j(a,b) 
=
\sum_{(\bfx,\bfy)\in\dcal_j(a,b)} \binom{j}{\bfx,\bfy} \prod_{i=1}^\alpha {\hat w}_{e_i}^{(x_i)}
\prod_{i=1}^\beta (-1)^{y_i }{\hat w}_{f_i}^{(y_i)}.
\end{equation*} 
Note that $\widehat{\fcal}_j(a,b)$ has at most $\depth^j$ terms (including the multinomial factor). Therefore,
since the function
\begin{equation*}
h(\mathbf{x}) = \prod_{j=1}^J x_j,
\end{equation*}
is continuously differentiable with bounded derivatives in $[-M^J, M^J]$, there is $C$
(depending on $M, J$) such that
\begin{equation*}
\sigma_j
\leq \tau^* + C \depth^j(2 \Sigma_{j-1}),
\end{equation*}
for small $\Sigma_{j-1}$.
Then we have
\begin{equation*}
\Sigma_J \leq \tau^* C^* \depth^{J^2/2},
\end{equation*}
for some $C^* > 0$ depending on $J,M$,
where we used $\sigma_2 \leq \tau^*$.

So it suffices to have $\tau^* = (\omega_n \depth^{J^2/2})^{-1}$ where
$\omega_n \to +\infty$ as $n\to +\infty$ arbitrarily slowly.
By the last part of Proposition~\ref{prop:concineq}, using a union
bound over the $O(n^2)$ short pairs of leaves, it follows that
$k = C'  \omega_n\depth^{2 J^2} \log n$ samples are enough to
guarantee
\begin{equation*}\label{eq:allmoments} 
\prob\left[
\left|
\delhat^{(j)} - \del^{(j)}
\right| \leq
(\omega_n \depth^{J^2/2})^{-1},\ \forall 1\leq j \leq J,\ \forall \mbox{ short pairs } a,b 
\right] 
\geq 1- o(1),
\end{equation*}
for some $C'$ depending on $J, M$.

As for the computational complexity of the algorithm, assume first that the tree is represented in such a way 
that finding the set of edges on the path between two leaves $a,b$ at distance $O(\depth)$ takes time
$O(\depth)$ (this is easy in a rooted tree). Note that for each $j, a ,b$ the sum
\begin{equation*}
\widehat{\fcal}_j(a,b) 
=
\sum_{(\bfx,\bfy)\in\dcal_j(a,b)} \binom{j}{\bfx,\bfy} \prod_{i=1}^\alpha {\hat w}_{e_i}^{(x_i)}
\prod_{i=1}^\beta (-1)^{y_i }{\hat w}_{f_i}^{(y_i)}.
\end{equation*} 
can be computed in time $\depth^J$. Since there are $O(n^2)$ pairs of leaves, 
the total complexity is $O(\depth^J n^2)$.
\end{proof}

Similarly, in the general case, we get:
\begin{proposition}\label{prop:concineqgen} 
Let $a, b, c\in L$ 
at graph distance less than $2 \depth$ where $\depth = \depth(T)$ is the depth
of $T$. Fix $j \in \nintgr$.
We have the following (where the constants depend on $J$ and $M$ only):
\begin{enumerate}
\item
There exists a constant $C$ such that, $\forall \lambda > 0$,
\begin{equation}\label{eqn:phihat} 
\prob\left(\left|{\hat\phi}_{ab|c}^{(j)}-\expec\left[{\hat\phi}_{ab|c}^{(j)}\right]\right| >
\lambda\right) \leq 2 \exp\left(-\frac{{\lambda}^2 k}{C \depth^{2j-1}}\right).
\end{equation}


\item There exists a constant $C'$ such that, if $k \geq \depth^2$,
\begin{equation}
\left|\expec\left[{\hat\phi}_{ab|c}^{(j)}\right] - \phi_{ab|c}^{(j)}\right| \leq C' \frac{(M\depth)^j}{\sqrt{k}}.
\end{equation}

\item If further
\begin{equation*}
C' \frac{(M\depth)^j}{\sqrt{k}} \leq \lambda,
\end{equation*}
then we have
\begin{equation}
\prob\left(\left|{\hat\phi}_{ab|c}^{(j)}-\phi_{ab|c}^{(j)}\right| >
2\lambda\right) \leq 2 \exp\left(-\frac{{\lambda}^2 k}{C \depth^{2j-1}}\right).
\end{equation}
\end{enumerate}
\end{proposition}
{\bf Proof Sketch:} The proof is very similar to Proposition~\ref{prop:concineq}. We only give a sketch.

To prove 1., it is enough to consider four separate cases depending on which path segment 
(corresponding to $H_1$, $H_2$, $H_3$ and $H_4$ in Figure~\ref{fig:oddrecursion}) we make the perturbation.
     
To prove 2., note that we can write
\begin{equation}\label{phi:alt-exp}
{\hat \phi}_{ab|c} = \frac{1}{k} \sum_1^k (X_i +\epsilon_I)^{j-1}(Y_i + \epsilon_{II}),
\end{equation}
with $X_i = (D_a^i -D_b^i) - \delta^{(1)}_{ab}$,
$Y_i = (D_a^i -D_c^i-\delta^{(1)}_{ac}) +(D_b^i-D_c^i-\delta^{(1)}_{bc})$, and 
\begin{eqnarray*}
\epsilon_{I} &=& \delta^{(1)}_{ab} - {\hat \delta}^{(1)}_{ab},\\
\epsilon_{II} &=& ( \delta^{(1)}_{ac} - {\hat \delta}^{(1)}_{ac}) +( \delta^{(1)}_{bc} - {\hat \delta}^{(1)}_{bc}).
\end{eqnarray*}
Also note that $\phi^{(j)}_{ab|c} = \expec[X_i^{j-1} Y_i], \forall i$. 
Use the Binomial theorem to expand the expression in (\ref{phi:alt-exp}) and write it as
\begin{equation*}
{\hat\phi}^{(j)}_{ab|c} = \frac{1}{k} \sum_{i=1}^k X_i^{j-1}Y_i + \rcal,
\end{equation*}
where the error term is
\[\rcal =  \epsilon_{II}\frac{1}{k}\sum_{i=1}^k (X_i +\epsilon_{I})^{j-1} + \frac{1}{k}\sum_{i=1}^k Y_i \sum_{l=1}^{j-1} \binom{j-1}{l} {\epsilon_I}^l X_i^{j-1-l}. \]
Now use the fact that $|X_i| \leq 4M \depth$, $|Y_i| \leq 8M\depth$, and Part 2.~of 
Proposition~\ref{prop:concineq} to conclude that 
\[\expec[|\rcal|] \leq C' \frac{(M\depth)^j}{\sqrt{k}}. \]
     
Part 3.~now follows by combining Part 1.~and 2.$\blacksquare$
\begin{theorem}\label{thm:generaler} 
Let $\eps > 0$ be arbitrarily small. If 
$k = \omega(\depth^{2J^2} \log n)$, then after an application of \textsc{ER}, one has
\begin{equation}
\prob \left[ \left|{\hat\delta}_e^{(j)}-\delta_e^{(j)}\right| 
\leq \eps,\ 
\forall e\in E,\ \forall 1\leq j \leq J \right] \geq 1-o(1),
\end{equation}
as $n \to +\infty$. The algorithm runs in time $O(\depth^{J} n^2)$.
\end{theorem}
\begin{proof}
The proof is identical to Theorem~\ref{thm:symer}.
\end{proof}

\section{Concluding Remarks}\label{sec:conclusion}

\begin{enumerate}

\item We have assumed that delays are {\em finitely supported}. 
This assumption is not essential. Unbounded distributions
for which similar concentration inequalities can be obtained
lead to the same results. For example, using~\cite[Proposition 4.18]{Ledoux:01},
one can treat the case of Exponential and Gamma delays.

\item It is an interesting problem, from a practical point of view, 
to improve the dependence of our results
on $J$.

\item It is somewhat intriguing that the reconstruction of the topology
of the tree required the joint distributions on {\em pairs} of leaves
whereas the reconstruction of delays (in the asymmetric case) required
the joint distributions on {\em triples} of leaves. A similar
situation holds in phylogenetics~\cite{Chang:96}. It could
be interesting to prove that this is indeed necessary in some sense.

\item Throughout, the model was assumed to be static. In real-life networks,
characteristics of the network change over time. One could try to adapt
our algorithm to a more dynamic setting. See for example~\cite{CaDaVaYu:00}
for a discussion of temporal issues.

\end{enumerate}


\section*{Acknowledgments}

We thank Gang Liang, Elchanan Mossel, and Bin Yu for discussions and encouragements. 
S.R. gratefully acknowledges the partial support of 
CIPRES (NSF ITR grant \# NSF EF 03-31494), NSERC (Canada), FQRNT (Quebec, Canada), a Lo\`eve Fellowship 
(UC Berkeley),
and NSF grant DMS-0528488. S.R. also thanks Martin Nowak and the
Program for Evolutionary Dynamics at Harvard University, where part of this work was done. 
S.B. thanks David Aldous for interesting discussions.

\bibliographystyle{amsplain}
\bibliography{all}

\providecommand{\bysame}{\leavevmode\hbox to3em{\hrulefill}\thinspace}
\providecommand{\MR}{\relax\ifhmode\unskip\space\fi MR }
\providecommand{\MRhref}[2]{%
  \href{http://www.ams.org/mathscinet-getitem?mr=#1}{#2}
}
\providecommand{\href}[2]{#2}
\begin{thebibliography}{10}

\bibitem{Atteson:99}
K.~Atteson, \emph{The performance of neighbor-joining methods of phylogenetic
  reconstruction}, Algorithmica \textbf{25} (1999), no.~2-3, 251--278.
  \MR{MR1703580 (2000k:92013)}

\bibitem{BhRaRo:arxiv}
S.~Bhamidi, R.~Rajagopal, and S.~Roch, \emph{Network delay inference from
  additive metrics}, Preprint. Available at Arxiv: math.PR/0604367, 2006.

\bibitem{Buneman:71}
P.~Buneman, \emph{The recovery of trees from measures of dissimilarity},
  Mathematics in the Archaelogical and Historical Sciences, Edinburgh
  University Press, Edinburgh, 1971, pp.~187--395.

\bibitem{CaDuHoTo:99}
Ram{\'o}n C{\'a}ceres, N.~G. Duffield, Joseph Horowitz, and Donald~F. Towsley,
  \emph{Multicast-based inference of network-internal loss characteristics},
  IEEE Trans. Inform. Theory \textbf{45} (1999), no.~7, 2462--2480.
  \MR{MR1725131}

\bibitem{CaDaVaYu:00}
Jin Cao, Drew Davis, Scott Vander~Wiel, and Bin Yu, \emph{Time-varying network
  tomography: router link data}, J. Amer. Statist. Assoc. \textbf{95} (2000),
  no.~452, 1063--1075. \MR{MR1821715}

\bibitem{CCLNY:04}
Rui Castro, Mark Coates, Gang Liang, Robert Nowak, and Bin Yu, \emph{Network
  tomography: recent developments}, Statist. Sci. \textbf{19} (2004), no.~3,
  499--517. \MR{MR2185628}

\bibitem{Chang:96}
Joseph~T. Chang, \emph{Full reconstruction of {M}arkov models on evolutionary
  trees: identifiability and consistency}, Math. Biosci. \textbf{137} (1996),
  no.~1, 51--73. \MR{MR1410044 (97k:92011)}

\bibitem{DaMoRo:09}
Constantinos Daskalakis, Elchanan Mossel, and S{\'e}bastien Roch,
  \emph{Phylogenies without branch bounds: Contracting the short, pruning the
  deep}, To appear in RECOMB'09. Preprint available as arXiv:0801.4190v1, 2009.

\bibitem{ErStSzWa:99a}
P{\'e}ter~L. Erd{\H{o}}s, Michael~A. Steel, L{\'a}szl{\'o}~A. Sz{\'e}kely, and
  Tandy~J. Warnow, \emph{A few logs suffice to build (almost) all trees. {I}},
  Random Structures Algorithms \textbf{14} (1999), no.~2, 153--184.
  \MR{MR1667319 (2000b:92003)}

\bibitem{Farris:73}
J.~S. Farris, \emph{A probability model for inferring evolutionary trees},
  Syst. Zool. \textbf{22} (1973), no.~4, 250--256.

\bibitem{Felsenstein:04}
J.~Felsenstein, \emph{Inferring phylogenies}, Sinauer, New York, New York,
  2004.

\bibitem{KiZhZh:03}
Valerie King, Li~Zhang, and Yunhong Zhou, \emph{On the complexity of
  distance-based evolutionary tree reconstruction}, Proceedings of the
  {F}ourteenth {A}nnual {ACM}-{SIAM} {S}ymposium on {D}iscrete {A}lgorithms
  ({B}altimore, {MD}, 2003) (New York), ACM, 2003, pp.~444--453. \MR{MR1974948}

\bibitem{LaceyChang:06}
Michelle~R. Lacey and Joseph~T. Chang, \emph{A signal-to-noise analysis of
  phylogeny estimation by neighbor-joining: insufficiency of polynomial length
  sequences}, Math. Biosci. \textbf{199} (2006), no.~2, 188--215. \MR{MR2211625
  (2007a:92048)}

\bibitem{Ledoux:01}
Michel Ledoux, \emph{The concentration of measure phenomenon}, Mathematical
  Surveys and Monographs, vol.~89, American Mathematical Society, Providence,
  RI, 2001. \MR{MR1849347 (2003k:28019)}

\bibitem{LiMoYu:u}
G.~Liang, E.~Mossel, and B.~Yu, \emph{Network topology inference through
  end-to-end measurements}, 2007.

\bibitem{Mossel:07}
E.~Mossel, \emph{Distorted metrics on trees and phylogenetic forests}, IEEE/ACM
  Trans. Comput. Bio. Bioinform. \textbf{4} (2007), no.~1, 108--116.

\bibitem{MosselRoch:06}
Elchanan Mossel and S{\'e}bastien Roch, \emph{Learning nonsingular phylogenies
  and hidden {M}arkov models}, Ann. Appl. Probab. \textbf{16} (2006), no.~2,
  583--614. \MR{MR2244426}

\bibitem{MotwaniRaghavan:95}
Rajeev Motwani and Prabhakar Raghavan, \emph{Randomized algorithms}, Cambridge
  University Press, Cambridge, 1995. \MR{MR1344451 (96i:65003)}

\bibitem{NiTatikonda:06}
J.~Ni and S.~Tatikonda, \emph{A {M}arkov random field approach to
  multicast-based network inference problems}, Proceedings of the IEEE
  International Symposium on Information Theory, 2006, pp.~2769--2773.

\bibitem{NiTatikonda:07}
\bysame, \emph{Explicit link parameter estimators based on end-to-end
  measurements}, Forty-Fifth Annual Allerton Conference, 2007.

\bibitem{NiTatikonda:08}
\bysame, \emph{Network tomography based on additive metrics}, Proceedings of
  the 42nd Annual Conference on Information Sciences and Systems, 2008,
  pp.~1149--1154.

\bibitem{NXTY:08}
Jian Ni, Haiyong Xie, S.~Tatikonda, and Y.R. Yang, \emph{Network routing
  topology inference from end-to-end measurements}, INFOCOM 2008. The 27th
  Conference on Computer Communications. IEEE (2008), 36--40.

\bibitem{LoDuHoTo:02}
Francesco~Lo Presti, N.~G. Duffield, Joe Horowitz, and Don Towsley,
  \emph{Multicast-based inference of network-internal delay distributions},
  IEEE/ACM Trans. Netw. \textbf{10} (2002), no.~6, 761--775.

\bibitem{SaitouNei:87}
N.~Saitou and M.~Nei, \emph{The neighbor-joining method: A new method for
  reconstructing phylogenetic trees}, Mol. Biol. Evol. \textbf{4} (1987),
  no.~4, 406--425.

\bibitem{SempleSteel:03}
Charles Semple and Mike Steel, \emph{Phylogenetics}, Oxford Lecture Series in
  Mathematics and its Applications, vol.~24, Oxford University Press, Oxford,
  2003. \MR{MR2060009 (2005g:92024)}

\bibitem{SneathSokal:73}
Peter H.~A. Sneath and Robert~R. Sokal, \emph{Numerical taxonomy}, W. H.
  Freeman and Co., San Francisco, Calif., 1973, The principles and practice of
  numerical classification, A Series of Books in Biology. \MR{MR0456594 (56
  \#14818)}

\bibitem{Vardi:96}
Y.~Vardi, \emph{Network tomography: estimating source-destination traffic
  intensities from link data}, J. Amer. Statist. Assoc. \textbf{91} (1996),
  no.~433, 365--377. \MR{MR1394093 (97a:62050)}

\end{thebibliography}

\clearpage

\appendix

\section{Examples of Regular Delay Distributions}\label{sec:regular}

Below, we give two typical examples of families of distributions
covered by our results. The first example is a set of continuous distributions
with few parameters. The second example is a general discrete distribution.
The latter is the main focus of~\cite{LoDuHoTo:02}.

\paragraph{Uniform distributions.} Let $\qcal = \{Q_\theta\}_{\theta\in \Theta}$
be the family of distributions where $Q_\theta$ is uniform on
$[0,\theta]$ with $\Theta = [\undtheta,\ovtheta]$ for some 
$0 < \undtheta < \ovtheta < +\infty$. Let ${\hat w}^{(2)}$ be the estimated
variance and define
\begin{displaymath}
{\hat\theta}^2
=
\Psi({\hat w}^{(2)})
=
\left\{
\begin{array}{ll}
\undtheta^2, & \mathrm{if\ }12 {\hat w}^{(2)} < \undtheta^2,\\
\ovtheta^2, & \mathrm{if\ }12 {\hat w}^{(2)} > \ovtheta^2,\\
12 {\hat w}^{(2)}, & \mathrm{otherwise.}
\end{array}
\right.
\end{displaymath}
Assume $|{\hat w}^{(2)} - w^{(2)}| \leq \delta \equiv \frac{\eps \undtheta^2}{12}$.
From $\theta^2 - {\hat\theta}^2 = (\theta - \hat\theta)(\theta + \hat\theta)$, it follows
easily that $|\theta - \hat\theta| \leq \frac{\eps \undtheta}{2}$. Note that
\begin{equation*}
\|Q_\theta - Q_{\hat\theta}\|_1
= \int_{0}^{\ovtheta} \left|\frac{\mathbf{1}_{x\leq \theta}}{\theta} 
- \frac{\mathbf{1}_{x\leq \hat\theta}}{\hat\theta}\right|\diff x,
\end{equation*}
and assuming w.l.o.g. that $\theta > \hat\theta$ (the other case is symmetric)
\begin{eqnarray*}
\int_{0}^{\ovtheta} \left|\frac{\mathbf{1}_{x\leq \theta}}{\theta} 
- \frac{\mathbf{1}_{x\leq \hat\theta}}{\hat\theta}\right|\diff x
= \hat\theta\left(\frac{1}{\hat\theta} - \frac{1}{\theta}\right) 
+ (\theta - \hat\theta)\frac{1}{\theta}
\leq 2\frac{\theta - \hat\theta}{\undtheta} \leq \eps.
\end{eqnarray*}
Therefore, $\qcal$ is $(\eps, 2)$-regular for any $\eps > 0$.

\paragraph{Bounded discrete distributions.}
Let $M$ be a positive integer and let $[M] = \{0,1,\ldots,M\}$. 
Also, let $0 < \undtheta < 1$ and
\begin{equation*}
\Theta = \left\{\theta = (\theta_0,\theta_1,\ldots,\theta_M)
: 0\leq \theta_i \leq 1,\ \forall i \in [M],
\ \theta_0 > \undtheta,\ \mathrm{and}\ \sum_{i\in[M]}i \theta_i \in [M] \right\}.
\end{equation*}
Denote by $\qcal = \{Q_\theta\}_{\theta\in \Theta}$
the family of distributions on $[M]$ such that
$X\sim Q_\theta$ means
\begin{equation*}
\prob[X = i] = \theta_i,\quad \forall i \in [M].
\end{equation*}
The assumption on the mean of $X$ in the definition of $\Theta$ greatly simplifies
the calculations below. It is a reasonable approximation in the standard practical case
where $\qcal$ is a discretization of continuous densities with a large
number of bins $M$. The assumption on $\theta_0$ simply indicates that the distribution
has been translated to ``start at 0.''
Define $\mu = \expec[X]$ where $X\sim Q_\theta$ and let
$\theta' = (\theta_{-M}', \theta_{-M+1}',\ldots, \theta_{M}')$ where
$\theta'_{i-\mu} = \theta_{i}$ for all $i\in[M]$ and $0$ otherwise. Note that
the following holds
\begin{equation*}
\sum_{i = -M}^{M} i^j \theta_i' = w^{(j)}(\theta),\quad \forall j\in[2M+1],
\end{equation*}
or in matrix form $\Lambda \theta' = \mathbf{w}$.
From the Vandermonde structure of $\Lambda$ it follows easily
that $\det \Lambda \geq 1$, that is, $\Lambda^{-1}$ exists, and furthermore
$\|\Lambda^{-1}\|_1$ is a strictly positive constant depending
on $\undtheta, M$. Let $\mathbf{\hat w}$ be the estimate of $\mathbf{w}$
and let $\hat\theta' = \Lambda^{-1}\mathbf{\hat w}$. Then, it follows that for
any $\eps > 0$ there is $\delta > 0$ such that
\begin{equation*}
\|\hat\theta' - \theta'\|_1 \leq \|\Lambda^{-1}\|_1 \|\mathbf{w} - \mathbf{\hat w}\|_1 \leq \eps,
\end{equation*}
whenever $\|\mathbf{w} - \mathbf{\hat w}\|_\infty \leq \delta$. Assume further
that $\eps < \undtheta/2$, then we can recover an estimate $\hat\theta$ of $\theta$ from $\hat\theta'$
such that $\|\hat\theta - \theta\|_1 \leq \eps$.
Indeed, our assumptions above allow us to infer a distribution centered at 0 which
we then translate to start at 0. Therefore, $\qcal$ is $(\eps, 2M-1)$-regular.
Note that strictly speaking one should force all components of $\hat\theta$ to be in $[0,1]$ and
renormalize appropriately. Details are omitted.

\clearpage

\section{DMR Algorithm}\label{section:dmr}

\begin{figure}
\begin{center}
\includegraphics[width=9cm]{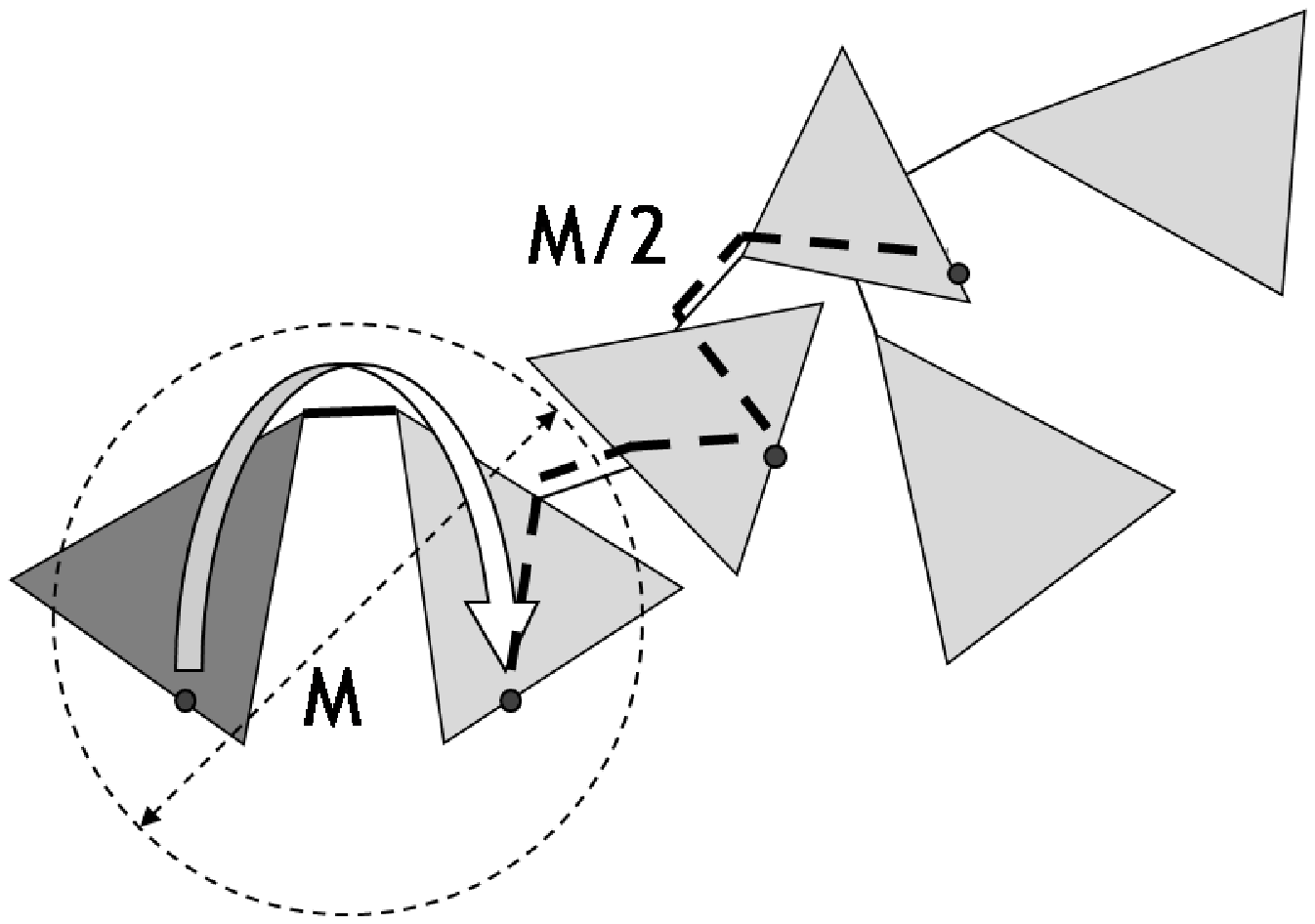}
\caption{Illustration of routine \extender.}\label{fig:extend}
\end{center}
\end{figure}

We shall now provide an outline of the DMR algorithm. 
The general DMR algorithm actually allows the user to build a ``forest'' when the number of
samples is too small. We will not use this feature here and we therefore simplify the algorithm
accordingly.
The input to the algorithm is a $(\tilde\tau, \widetilde M)$-distorted
metric $\distdw$ on $n$ leaves. 
In particular, we assume that the values $\tilde\tau$ and
$\widetilde M$ are known to the algorithm. 
We denote the true tree by $T = (V,E)$. 
Take $\alpha, \alpha' > 0$
and $0 < \beta, \beta' < 1$ such that
\begin{equation*}
6 < \alpha' + 3 < \alpha < (\tilde\alpha)^{-1},
\end{equation*}
and
\begin{equation*}
(\tilde\beta)^{-1} \widetilde M + \tilde\tau < \beta \widetilde M < \frac{1}{2}[\beta' \widetilde M - 3\tilde\tau].
\end{equation*}
(Here it is assumed that $\widetilde M = \omega(\tilde \tau)$.)
The details of the subroutines \minicontractor~and \extender~can be found in Figures~\ref{figure:minicontractor} 
and~\ref{figure:extender}.
The reader is referred to~\cite{DaMoRo:09} for a
detailed explanation of the algorithm---which is somewhat involved.
In a nutshell, for each pair of leaves $u,v$ that are not ``too far'':
1) the algorithm finds all edges sitting on the path between $u$ and $v$
(as illustrated in Figure~\ref{fig:mini});
2) then it derives the bipartitions corresponding to these edges
by ``extending'' the bipartitions in a small ball around $u,v$
(as illustrated in Figure~\ref{fig:extend}).
\begin{itemize}
\item \textbf{Pre-Processing: Proximity Test.} 
Build the graph
$\cluster{\beta} = (\vcluster{\beta}, \ecluster{\beta})$ where $\vcluster{\beta} = L$
and
$(u,v) \in \ecluster{\beta} \iff \distdw(u,v) < \beta \widetilde M$;
\item \textbf{Main Loop.} 
\begin{itemize}
\item For all pairs of leaves
$u,v \in \vcluster{\beta}$ such that $(u,v) \in \ecluster{\beta}$:  
\begin{itemize}
\item \textbf{Mini Reconstruction.}
Compute 
\begin{equation*}
\{\minipart{j}{u}{v}\}_{j=1}^{\partsize{u}{v}} := \mathrm{\minicontractor}(\cluster{\beta}; u,v);
\end{equation*}
\item \textbf{Bipartition Extension.}
Compute 
\begin{equation*}
\{\fullpart{j}{u}{v}\}_{j=1}^{\partsize{u}{v}} := \mathrm{\extender}(\cluster{\beta}, \{\minipart{j}{u}{v}\}_{j=1}^{\partsize{u}{v}}; u,v);
\end{equation*}
\end{itemize}
\item Deduce the tree $\widehat{T}$ from $\{\fullpart{j}{u}{v}\}_{j=1}^{\partsize{u}{v}}$;
\end{itemize}
\item \textbf{Output.} Return the resulting tree $\widehat{T}$.
\end{itemize}
\begin{figure}
\begin{center}
\input{mini-tomo.pstex_t}\caption{Illustration of routine \minicontractor.}\label{fig:mini}
\end{center}
\end{figure}

\begin{figure*}
\framebox{
\begin{minipage}{16.5cm}
{\small \textbf{Algorithm} \minicontractor\\
\textit{Input:} Graph $\cluster{\beta}$; Leaves $u,v$;\\
\textit{Output:} Bipartitions $\{\minipart{j}{u}{v}\}_{j=1}^{\partsize{u}{v}}$;
\begin{itemize}
\item \textbf{Ball.} Let
\begin{equation*}
\ball{\beta'}{0}(u,v) := \left\{w \in \cluster{\beta}\ :\ 
\distdw(u,w) \lor \distdw(v,w) < \beta' \widetilde M \right\};
\end{equation*}
\item \textbf{Intersection Points.} For all $w\in \ball{\beta'}{0}(u,v)$, estimate the point of intersection
between $u,v,w$ (distance from $u$), that is,
\begin{equation*}
\intersect_w := \frac{1}{2}\left(
\distd(u,v) + \distd(u,w) - \distd(v,w)
\right);
\end{equation*}
\item \textbf{Long Edges.} 
Set $S := \ball{\beta'}{0}(u,v) - \{u\}$, $x_{-1} = u$, $j := 0$; 
\begin{itemize}
\item Until $S=\emptyset$:
\begin{itemize}
\item Let $x_0 = \arg\min\{\intersect_w\ :\ w \in S\}$ (break ties arbitrarily); 
\item If $\intersect_{x_0} - \intersect_{x_{-1}} \geq \alpha' \tilde\tau$, create a new edge
by setting $\minipart{j+1}{u}{v} := \{\ball{\beta'}{0}(u,v) - S, S\}$ and let $C_{j+1} := \{x_0\}$, $j := j+1$; 
\item Else, set $C_j := C_j \cup \{x_0\}$;
\item Set $S := S-\{x_0\}$, $x_{-1} := x_0$;
\end{itemize}
\end{itemize}
\item \textbf{Output.} Return the bipartitions $\{\minipart{j}{u}{v}\}_{j=1}^{\partsize{u}{v}}$.
\end{itemize}
}
\end{minipage}
} \caption{Algorithm \minicontractor.} \label{figure:minicontractor}
\end{figure*}

\begin{figure*}
\framebox{
\begin{minipage}{16.5cm}
{\small \textbf{Algorithm} \extender\\
\textit{Input:} Graph $\cluster{\beta}$; 
Bipartitions $\{\minipart{j}{u}{v}\}_{j=1}^{\partsize{u}{v}}$; Leaves $u,v$;\\
\textit{Output:} Bipartitions $\{\fullpart{j}{u}{v}\}_{j=1}^{\partsize{u}{v}}$;
\begin{itemize}
\item For $j=1,\ldots,\partsize{u}{v}$ (unless $\partsize{u}{v}=0$):
\begin{itemize}
\item {\bf Initialization.} Denote by $\minipartone{j}{u}{v}$ the vertex set containing $u$ in the bipartition $\minipart{j}{u}{v}$,
and similarly for $v$; Initialize the extended partition $\fullpartone{j}{u}{v} := \minipartone{j}{u}{v}$,
$\fullparttwo{j}{u}{v} := \miniparttwo{j}{u}{v}$;
\item {\bf Modified Graph.} Let $K$ be $\cluster{\beta}$ where all edges between
$\minipartone{j}{u}{v}$ and $\miniparttwo{j}{u}{v}$ have been removed;
\item {\bf Extension.} For all $w \in \vclusters{\beta}{i} - (\minipartone{j}{u}{v}\cup\miniparttwo{j}{u}{v})$,
add $w$ to the side of the partition it is connected to in
$K$ (by definition of $K$, each $w$ as above is connected to exactly
one side); 
\end{itemize}
\item Return the bipartitions $\{\fullpart{j}{u}{v}\}_{j=1}^{\partsize{u}{v}}$.
\end{itemize}
}
\end{minipage}
} \caption{Algorithm \extender.} \label{figure:extender}
\end{figure*}

\end{document}